\documentclass[oneside,reqno,a4paper,11pt]{amsart}
\usepackage{amssymb}
\usepackage{hyperref}

\newtheorem {theorem}{Theorem}[section]
\newtheorem {lemma}[theorem]{Lemma}

\newtheorem {proposition}[theorem]{Proposition}
\theoremstyle{remark}
\newtheorem {remark}{Remark}[section]

\theoremstyle{problem}

\theoremstyle{plain} \numberwithin {equation}{section}

\everymath{\displaystyle}

\allowdisplaybreaks

\makeindex

\def\Xint #1{\mathchoice
{\XXint \displaystyle \textstyle {#1}} %
{\XXint \textstyle \scriptstyle {#1}} %
{\XXint \scriptstyle \scriptscriptstyle {#1}} %
{\XXint \scriptscriptstyle \scriptscriptstyle {#1}} %
\!\int}
\def\XXint #1#2#3{{\setbox 0=\hbox {$#1{#2#3}{\int }$}
\vcenter {\hbox {$#2#3$}}\kern -.5\wd 0}}

\def\dashint {\Xint -}
\def\be{\begin{equation}}
\def\ee{\end{equation}}
\def\ba{\begin{aligned}}
\def\ea{\end{aligned}}

\def\g{\nabla}
\def\lap{\Delta}
\def\O{\Omega}
\def\o{\omega}

\def\R{\mathbb{R}}
\def\p{\partial}
\def\a{\alpha}

\def\f{\frac}
\def\th{\theta}
\def\d{\delta}

\def\th{\theta}
\def\div{\mathrm{div}}
\def\ld{\lambda}

\def\endproof{{\hfill$\Box$}\\}

\def\eps{\varepsilon}

\def\trho{\overline{\rho}}
\def\mrho{{\widehat{\rho}}}
\def\tp{\overline{p}}
\def\tC{\Lambda}
\def\t#1{\overline{#1}}
\def\into{\int_{\O}}
\def\aint{\dashint_{\O}}

\def\E{{\mathcal{E}}}

\def\defeq{\triangleq}
\begin{document}

 %{\noindent\small Draft}
\vspace{1cm}
\title[Global weak solutuons to Dirichlet problems of CNS]
{Global weak solutions with higher regularity to the compressible Navier-Stokes equations under Dirichlet boundary conditions}

\author{Xiangdi, Huang}
\author{Zhouping Xin}
\author{Wei YAN}

\begin{abstract}
In this manuscript, we aim to establish global existence of weak solutions with higher regularity to the compressible Navier-Stokes equations under no-slip boundary conditions. Though Lions\cite{L1} and Feireisl\cite{F1} have established global weak solutions with finite energy under Dirichelet boundary conditions by making use of so called effective viscous flux and oscillation defect measure,Hoff has  investigated global weak solutions with higher regularity in \cite{H1,Hof2} when the domain is either whole space or half space with Navier-slip boundary conditions, yet the existence theory of global weak solution with higher regularity under Dirichlet boundary conditions remains unknown. In this paper we prove that the system will admit at least one global weak solutions with higher regularity as long as the initial energy is suitably small when the domain is a 2D solid disc. This is achieved by exploiting the structure of the exact Green function of the disc to decompose the effective viscous flux into three parts, which corresponds to the pressure term, boundary term and the remaining term respectively. In order to control the boundary term, one of the key observations is to use the geometry of the domain which sucessfully to bound the integral of the effective viscous flux where $L^1$ norm is always unbounded.

\end{abstract}

\address{Xiangdi Huang\hfill\break\indent
Institute of Mathematics, \hfill\break\indent
Academy of Mathematics and systems science, Chinese Academy of Sciences, Beijing 100190, P. R. China}
\email{\href{xdhuang@amss.ac.cn}{xdhuang@amss.ac.cn}}

\address{Zhouping Xin\hfill\break\indent
The Institute of Mathematical Sciences,\hfill\break\indent
The Chinese University of Hong Kong, Shatin, N.T., Hong Kong.}
\email{\href{zpxin@ims.cuhk.edu.hk}{zpxin@ims.cuhk.edu.hk}}

\address{Wei YAN\hfill\break\indent
School of Mathematics,\hfill\break\indent
Jilin University, Changchun, P.R.China}
\email{\href{wyanmath@gmail.com}{wyanmath@gmail.com}}
\maketitle
\def\q0{{q_0}}
\section{Introduction}

The study of the isentropic compressible Navier-Stokes equations with constant viscosity coefficients has a long history. There are plenty of literatures concerning the existence theory of local strong and classical solutions to the multi-dimensional compressible Navier-Stokes equations where the density has a positive lower bound, see the references \cite{Be, Itaya, Ka-1, Ka-2, Na, Salvi, Solo, Valli-1, Valli-2} and therein. The local well-posedness of strong and classical solutions allowing initial vacuum has been obtained by Cho and Kim in \cite{K1,K3,K2}.

 In their recent paper, Merle\cite{FM} et al.  proved that the smooth solution of the 2D radially symmetric compressible Navier-Stokes equation will inevitably form a shell singularity in finite time for non-vacuum data. It is also shown by Xin\cite{X1} and Xin-Yan\cite{XY} that the classical solution does exist globally in general in the presence of vacuum. In this respect, in order to study the long time dynamics of the fluid and capture those solutions with lower regularities, the investigation of weak solutions becomes a key issue.

 A major breakthrough is due to Lions\cite{L1}. By fully making use of find properties of so called effective viscous flux, he proved the global existence of weak solutions allowing initial vacuum of the three dimensional isentropic compressible Navier-Stokes equations with "finite energy" as long as the adiabatic exponent $\gamma> \frac{9}{5}$.Later, Fereisl\cite{F1} introduced oscillation defect measure to study renormalized weak solutions and succeeded to  improve the restrcitions on the adiabatic exponent to $\gamma>\frac{3}{2}$. Whether the compressible Navier-Stokes equations admit global weak solutions under the physical condition $\gamma>1$ becomes an outstanding chanllenging problem. This was only achieved for the spherically symmetric flows with large initial data by Jiang-Zhang\cite{JZ} when the domain is the whole space $R^3$, or the solution with "small initial energy" by Hoff\cite{H1,Hof3,Hof2} for either whole space of half space with Navier-slip boundary conditions. However, the uniqueness of weak solutions remains unknown.

Indeed,  the regularity of weak solutions constructed by Hoff\cite{H1} is stronger than Lions-Feireisl's solutions in the sense that particle path can be defined in the non-vacuum region among others and weaker than the usual strong solutions in the sense that discontinuity of the density can be transported along the particle path. Unfortunately, it still remains unkown that whether Lions-Feireisl's weak solution is unique.  The reason why we are interested in looking for weak solutions with higher regularity is related to an chanllenging problem, which may guarantee the uniqueness. For example, this was verified for the piecewise smooth solution\cite{Hof3} and the spherically symmetric solution\cite{Hxd-sy}.

Motivated by the same purpose, Desjardins\cite{Des} proved the local weak solutions with the same regularities of Hoff\cite{H1}'s small energy solution for the periodic domain $T^3$ by asking $\gamma>3$. It was extended to the physical restriction for $\gamma>1$ by Huang-Yan\cite{Hy} for periodic domains. In both cases, the weak solutions $(\rho, u)$ to the compressible isentropic Navier-Stokes equations enjoy more regularities as follows
\be\label{reg}
\|\rho\|_{L^\infty(0,T;L^\infty(\Omega))} + \|\nabla u\|_{L^\infty(0,T; L^2(\Omega))}\le C.
\ee

However, the existence theory of global weak solution with higher regularity \eqref{reg} under Dirichlet boundary conditions remains unknown. It is easy to see that the effective viscous flux satifies an elliptic equations. One can derive global estimates over the whole domains when the domain is whole space, periodic or some bounded one with Navier-slip boundary conditions. But unfortunately, this becomes a big issue for the important no-slip boundary conditions. Namely, the lack of boundary condition of the effective viscous flux will is the main obstacles preventing one obtaining upper bound of the density in this case. This is the main problem we will address here.

A general viscous isentropic compressible fluid occupying a domain $\O\subset R^N$ is governed by the compressible Navier-Stokes equations as follows
\be\label{ns-0}
\left\{\ba
&\p_t\rho + \div(\rho u) = 0,\\
&(\rho u)_t + \div(\rho u\otimes u) +\g p = \mu\lap u + \ld(\g\div u),
\ea\right.\ee
where $\rho\ge 0$, $u=(u^1,u^2,...u^N)$, and $p=a\rho^\gamma(a>0,\gamma>1)$ are the fluid density, velocity, and pressure, respectively. The viscosity coefficients $\mu$ and $\lambda$ are assumed to be constants satisfying
\be
\mu>0,\quad \mu+\frac{N}{2}\lambda\ge 0.
\ee
In this paper, $\Omega$ will be taken to  be a two-dimensional unit ball in $R^2$.

We are looking for the solutions $(\rho(x,t),u(x,t))$ to the initial-boundary value problem for the system \eqref{ns-0}  with the initial-boundary data
\be
 (\rho,u)|_{t=0}= (\rho_0,u_0),\quad  u\big|_{\p\O} = 0.
\ee

First we will list the notations and conventions used frequently  in this paper.

The average of a function $\bar{f}$ is defined as
\be
\bar{f} \defeq\aint f(x)dx.
\ee
Thus the average of density and pressure $\bar{\rho}(t)$ and  $\bar{p}(t)$ can be defined as follows
\be
\trho(t) \defeq \aint \rho(x,t)dx,\quad \tp(t) \defeq \aint p(x,t)dx.
\ee
By the conservation of mass, $\trho(t)=\trho$ is indeed a constant corresponding to the initial total mass.

The initial energy $E_0$ is defined as
\be
E_0 \defeq \into\f 12 \rho_0|u_0|^2 + G(\rho_0) dx.
\ee
where G denotes the potential energy density given by
\be
G(\rho) \defeq \rho\int_{\trho}^\rho \f{p(s)-p(\trho)}{s^2}ds=\f{\rho^\gamma - \rho\trho^{\gamma-1}}{\gamma-1} + (\trho-\rho)\trho^{\gamma-1}.
\ee
It is clear that
\be
 \left\{\ba
 &G(\rho )=\frac{1}{\gamma-1}p, \quad\quad\quad\quad\text{for}\  \bar\rho=0,\\
&G(\rho )\geq C(\mrho,\trho) |\rho-\trho|^2,\quad\text{for } 0\leq \rho\leq\mrho,\\
&G(\rho )\geq C(\mrho,\trho) |\rho-\trho|^{\gamma},\quad\text{for } \mrho\leq\rho <\infty.
\ea\right.
\ee
Moreover, we use $C$ as a uniform constant which may varies in different places and, we may specify $C$ as $C_1,\ C_2$ through out this paper.

Multiplying the momentum equations $\eqref{ns-0}_2$ by $u$ and making use of the density equation $\eqref{ns-0}_1$, one can derive easily the following basic energy estimate of the system \eqref{ns-0}.

\be
 \into\left(\frac{1}{2} \rho|u|^{2}+G(\rho)\right) d x
+\int_{0}^{T} \into\left(\mu|\g u|^{2}+\ld(\div u)^{2}\right) d x d t = E_{0}.
\ee
In this paper, we will show that the system will admit at least one global weak solutions with higher regularity as long as the initial energy is suitably small when the domain is a 2D solid ball. We will use the exact form of the Green function of the ball to decompose the effective viscous flux into three parts, which corresponds to the pressure term, boundary term and the remaining term respectively. In order to control the boundary term, one of the key observations is to use the geometry of the domain to bound the integral of the effective viscous flux whose $L^1$ norm is always unbounded.

Before stating our main theorem, we first introduce the following norms and notations. For any function $f$ and a non-negative integer $k\ge 0,p\ge 1$, the $L^p$ and $W^{k,p}$ norm of the function $f$ is defined as follows
\be
\|f\|_{L^p(\Omega)}=\Big(\int_{\Omega} f^pdx\Big)^{\frac{1}{p}}, \|f\|_{W^{k,p}(\Omega)}=\sum_{i=0}^k\|D^if\|_{L^p}.
\ee
Here
\begin{equation*}
L^p=L^p(\Omega), \quad W^{k,p}=\{f\in L_{loc}^1(\Omega)|  \|\nabla^k f\|_{L^p}<\infty\}, H^k=W^{k,2}.
\end{equation*}
Morover, for any number $s\in(0,1)$ and $1\le p<\infty$, the function space $W^{s,p}(\Omega)$ is defined as
\be
W^{s,p}(\Omega)=\Big\{f\in L^p(\Omega): \int_{\Omega}\int_{\Omega}\frac{|f(x)-f(y)|^p}{|x-y|^{n+sp}}dxdy<\infty\Big\},
\ee
and
\be
H^s(\Omega)=W^{s,2}(\Omega).
\ee
The norm of a  function $f$ in $W^{s,p}(\Omega)(s\in(0,1),p\in[1,\infty))$ is defined by
\be
\|f\|_{W^{s,p}(\Omega)}^p = \|f\|_{L^p}^p + \int_{\Omega}\int_{\Omega}\frac{|f(x)-f(y)|^p}{|x-y|^{n+sp}}dxdy.
\ee

{\bf Assumptions:}  Assume that $\gamma\in (1,2)$. Denote $\q0 \defeq \f{12\gamma}{\gamma-1}$. For any given $\beta\in \big(\f 12, 1\big]$, $\d_0 \defeq \f{2\beta-1}{3\beta}\in\big(0, \f 13\big]$. Assume that $\mrho \geq \trho +1$, and
\be
\|\rho_0\|_{\q0}\leq\mrho,\quad\|u_0\|_{\dot H^{\beta}}\leq M.
\ee
Now the main results of this paper can be stated as follows.
\begin{theorem} \label{th-main}
Let $\Omega$ be a two-dimensional unit disc. Assume that the initial data $(\rho_0,u_0)$ satisfy the following regularity
\be
0\le \rho_0\in L^{q_0}(\O),\ u_0\in H^{\beta}(\O).
\ee
Assume that $\f {2\mu}{2\mu+\ld}\leq \min\left\{\f 1{4\tC(\f{\gamma+11}{\gamma-1})},\ \f 1{4 \tC^{\f 14}(4)}\right\}$,where $\tC(q)$ is the constant defined in Lemma \eqref{lem-Ag}.

Then there exists $\eps >0$, depending on $\mu, \ld, \gamma, \q0, \beta, \d_0, M, \trho, \mrho$ such that, if
\be
E_0\leq \eps.
\ee
the Dirichlet problem of the system \eqref{ns-0} admits at least one global weak solution satisfying
\be
%\left\{
\ba
&0\le\rho\in L^\infty([0,\infty);L^{q_0})\cap C([0,\infty);L^{q}),\quad \mbox{for all $q\in[1,q_0)$},\\
& u\in C([0,\infty); H^{\beta}),\quad (\rho\dot{u},\ \sigma(t)^{\frac{1}{2}}\nabla\dot{u})\in L^2((0,T)\times\Omega).
%\right.
\ea
\ee
where $\sigma(t)=\min(1,t)$.
\end{theorem}
\begin{remark}
The regularity of weak solution established here is between Lions-Feiresl's finite energy solution and Hoff's small energy solution.
\end{remark}
\begin{remark}
In fact, Pereplista\cite{per} first investigated the weak solution to the multi-dimensional compressible Navier-Stokes equations under Dirichlet boundary condition. When the density is away from vacuum and  $C^{\alpha}$ near boundary, he showed the global existence of weak solution. But his proof relies heavily on two key ingredients: one is the absence of vacuum, the other is the Holder continuity of the density in a neighborhood of the boundary, thus his method cannot be applied to our case. We need to carefully analysis the exact structure of the Green function of the domain, which is necessary to decompose the viscous effective flux into three parts. Based on this new observation, we can establish weighted estimate of the density.
\end{remark}

Denoting that
\be
A_{1}(T) \defeq \sup _{t \in[0, T]}\left(\sigma\|\g u\|_{2}^{2}\right)+\int_{0}^{T} \into \sigma \rho|\dot{u}|^{2} d x d t
\ee
and
\be
A_{2}(T) \defeq \sup _{t \in[0, T]} \sigma^{3} \into \rho|\dot{u}|^{2} d x+\int_{0}^{T} \into \sigma^{3}|\g \dot{u}|^{2} d x d t
\ee
and
\be
A_{3}(T) \defeq \sup _{0 \leq t \leq T} \into \rho|u|^{3}(x, t) d x.
\ee

The key idea in the proof of Theorem \ref{th-main} is the following bootstrap argument.
\begin{proposition}\label{prop-1}
Under the condition of Theorem 1.1, there exists some positive constant $\epsilon$ depending on $\mu, \ld, \gamma, \q0, \beta, \d_0, M, \trho, \mrho$ such that if $(\rho,u)$ is a solution of \eqref{ns-0} satifying
\be\label{H1}
\sup\limits_{0\leq t\leq T}\|\rho\|_{\q0}\leq 2\mrho,\quad A_1(T)+A_2(T)\leq 2E_0^{\f 1{2\gamma}},\quad
A_3(\sigma(T)) \leq 2E_0^{\d_0}.
\ee
Then the following estimates hold:
\be
\sup\limits_{0\leq t\leq T}\|\rho\|_{\q0}\leq \f 74 \mrho
\ee
and
\be
A_1(T)+A_2(T)\leq E_0^{\f 1{2\gamma}},\quad A_3(\sigma(T)) \leq E_0^{\d_0}.
\ee
provided  $E_0\leq \eps$.
\end{proposition}

\section{Preliminaries}
In this section, we list some basic ellliptic estimates and elementary decompositions of the effective viscous flux and the corresponding estimates.
\subsection{Decomposition of the velocity}

In the rest of the paper, we always using the notation $C(\mrho)=C(\mrho, \trho)$ for simplicity.

Define $v$ and $w$ as follows:
\be
\left\{\ba
&\mu\lap v +\ld\g\div v = \g p,\\
&v\big|_{\p\O}=0,
\ea\right.
\qquad
\left\{\ba
&\mu\lap w +\ld\g\div w = \rho \dot u,\\
&w\big|_{\p\O}=0,
\ea\right.
\ee
Introduce the Lame operator $L$ as $L f =\mu\lap f +\ld\g\div f$. And denote $L^{-1}$ the inverse operator of $L$ with Dirichlet boundary condition $f(\p\O)=0$.

\begin{lemma} \label{lem-du}
For any $1<q<\infty$, it holds that
\be\label{l-1}
\|\g v\|_q\leq C(q)\|p-\tp\|_q,\quad
\|\g^2 w\|_q\leq C(q)\|\rho \dot u\|_q,
\ee
and
\be
\|\g u\|_q \leq C(q) \big(\|\rho^{\f 12}\|_q\|\rho^{\f 12}\dot u\|_2  + \|p-\tp\|_q\big).
\ee
\end{lemma}
\proof
\eqref{l-1} follows immediately from the elliptic theory to the Lame operator. Moreover,
\[\ba
\|\g u\|_q &\leq \|\g w\|_q + \|\g v\|_q
\leq C(q)\|\g^2 w\|_{\f{2q}{q+2}} + C(q)\|p-\tp\|_q \\
&\leq C(q)\|\rho\dot u\|_{\f{2q}{q+2}} + C(q)\|p-\tp\|_q \\
&\leq C(q) \big(\|\rho^{\f 12}\|_q\|\rho^{\f 12}\dot u\|_2  + \|p-\tp\|_q\big)\\
\ea\]
\endproof

\subsection{Decomposition of the effective viscous flux}
\be
\quad
\ee
In this section, we will study the decomposition of the effective viscous flux in a general Lame system, i.e
\be\label{lame-1}
\left\{\ba
&\mu\lap u + \ld\g\div u = \g g + f,\\
& u=0,\quad\text{ on }\p\O.
\ea\right.\ee
The effective viscous flux F to the above system \eqref{lame-1} is defined as
\be\label{ff-0}
F \defeq (\mu + \ld)\div u - g.
\ee
Then F satisfies
\be\label{ff-1}
\lap F = \div f.
\ee

Assume that $\bar{G}(x,y)\in C^\infty(\Omega\times\Omega\mbox{$\backslash$}D)$ with $D=\{(x,y)\in \Omega\times\Omega |\ (x=y)\}$ is Green function on $\O$. Then a general function $H$ has the following explict expression
\be\label{ff-2}
\ba
H(x) &= \int_\O \bar{G}(x,y) \lap H(y) dy + \int_{\p\O} \f{\p \bar{G}(x,y)}{\p n} H(y) dy\\
& \defeq G (\lap H) + G_b H.
\ea\ee
where $G$ and $G_b$ are operators defined as
\[
G(H)(x) =\bar{G}*H= \int_\O \bar{G}(x,y) H(y) dy,\quad G_b H(x) = \int_{\p\O} \f{\p \bar{G}(x,y)}{\p n} H(y) dy.
\]
Assume that $u$ is a solution of the Lame system \eqref{lame-1}. Then it holds that
\be
u=G(\lap u).
\ee
Then in view of \eqref{ff-1}-\eqref{ff-2}
\be\label{ff-3}
F = G(\div f) + G_b(F),
\ee
and by definition \eqref{ff-0}
\be\ba
\f{\ld}{\mu + \ld}\g F &= \ld\g\div u - \f{\ld}{\mu + \ld}\g g = \g g + f - \mu\lap u - \f{\ld}{\mu + \ld}\g g\\
&=-\mu\lap u + \f{\mu}{\mu + \ld}\g g + f
\ea\ee
The following equality also holds true
\be
\lap(x F) = 2\g F + x\div f
\ee
Consequently,
\[\ba
\f{\ld}{2(\mu + \ld)}\lap(x F) &+ \mu\lap u = \f{\ld}{\mu + \ld}\g F + \f{\ld}{2(\mu + \ld)}x\div f + \mu\lap u\\
&= \f{\ld}{\mu + \ld}\g F + \f{\ld}{2(\mu + \ld)}x\div f + \mu\lap u\\
&=\f{\ld}{2(\mu + \ld)}x\div f + \f{\mu}{\mu + \ld}\g g + f\\
&\defeq Q.
\ea\]
Therefore the effective viscous flux F can be expressed as
\be\label{gg-1}
\f{\ld}{2(\mu + \ld)} x F + \mu u =G(Q) + \f{\ld}{2(\mu + \ld)} G_b(x F).
\ee
Taking $div$ on both sides of \eqref{gg-1} yields
\be
\f{\ld}{2(\mu + \ld)} (nF + x\cdot\g F) + \mu\div u =\div G(Q) + \f{\ld}{2(\mu + \ld)} \div G_b(x F).
\ee
By \eqref{ff-3}
\be
\g F = \g G(\div f) + \g G_b (F)
\ee
and
\be
\mu\div u = \f{\mu}{\mu+\ld}(F+g)
\ee
one gets
\be\ba
\f{n\ld}{2(\mu + \ld)} F &+ \f{\ld}{2(\mu + \ld)}x\cdot\g G(\div f) + \f{\ld}{2(\mu + \ld)}x\cdot\g G_b (F) + \f{\mu}{\mu+\ld}(F+g)\\
&=\div G(Q) + \f{\ld}{2(\mu + \ld)} \div G_b(x F).
\ea\ee
Therefore
\be\label{ff-4}
\ba
\f{n\ld+2\mu}{2(\mu + \ld)} F &+ \f{\mu}{\mu+\ld} g =\div G(Q) - \f{\ld}{2(\mu + \ld)}x\cdot\g G(\div f) \\
&+ \f{\ld}{2(\mu + \ld)} (\div G_b(x F) - x\cdot\g G_b (F)).
\ea\ee

Finally, we have the following decomposition of the viscous effective flux F.
\begin{lemma} \label{lem-Farb}
Assume that $\O$ is two-dimensional unit disc, $F$ is the effective viscous flux corresponding to the Lame system \eqref{lame-1}. Then F has the following decomposition
\be\label{F-1}
\ba
 F = \f{2\mu}{2\mu + \ld} A(g) + B(F) +R(f).
\ea\ee
with
\be\label{arb}\left\{\ba
A(g) &= \div G(\g g)-g,\\
B(F) &= \f{\ld}{2\pi R(2\mu + \ld)}\int_{\p\O} F(y) dy,\\
R(f) &= \f{2(\mu + \ld)}{2\mu+\ld}\div G(f) - \f{\ld}{2\mu + \ld}\Big(\div G(x\div f)-x\cdot\g G(\div f) + G(\div f)\Big).
\ea\right.\ee
where A is related to the pressure term, B is linked to the boundary integral of the  effective viscous flux, R is the remaining source term f. In the original equations \eqref{ns-0}, R is corresponding to the material derivatives $\rho\dot{u}$.
\end{lemma}
\begin{remark}
As it will be seen later, the most difficult part is how to treat the bound for A and B in the original Navier-Stokes equation. This is the key issue we will pay most attention to. The lack of boundary information of F is the crucial obstacle preventing one for deriving the upper bound of the density. But fortunately, by a delicate analysis of A and new observations on the boundary integral B due to the special geometry, we can obtain the $L^p$ bound of the density for suitably large $p\in(1,\infty)$.

\end{remark}
\proof
Note that the Green function $G(x,y)$ of the n-dimensional ball with radius R is given by
\be
G(x,y)=\left\{
\ba
&\frac{1}{(2-n)\omega_n}\big(|x-y|^{2-n}-\big|\frac{R}{|x|}x-\frac{|x|}{R}y \big|^{2-n}\big)\quad\mbox{for $n\ge 3$}, \\
&\frac{1}{2\pi}\big(\log|x-y|-\log\big|\frac{R}{|x|}x-\frac{|x|}{R}y\big|\big)\quad\mbox{for $n=2$}
\ea
\right.
\ee
where $\omega_n$ denotes the surface area of the unit sphere in $\R^n$.

It is easy to verify that for any $n\ge 2$
\[
\f{\p G(x, y)}{\p n}\Big|_{\p\O} = \f{R^2-|x|^2}{n\o_nR|x-y|^n}.
\]
and
\[
 \f{\p \g_x G(x, y)}{\p n}\Big|_{\p\O} = \f{-2x}{n\o_n R|x-y|^n} -  \f{R^2-|x|^2}{\o_n R|x-y|^{n+1}}\f{x-y}{|x-y|}
\]
So
\[\ba
(y-x)\cdot\f{\p \g_x G(x, y)}{\p n}\Big|_{\p\O} &= \f{-2x(y-x)}{n\o_n R|x-y|^n} -  \f{R^2-|x|^2}{\o_n R|x-y|^{n+1}}\f{(x-y)(y-x)}{|x-y|}\\
%&=\f{-2x(y-x)}{n\o_n R|x-y|^2} +  \f{R^2-|x|^2}{\o_n R|x-y|^n}=\f{nR^2-n|x|^2 -2x(y-x)}{n\o_n R|x-y|^n}\\
%&=\f{(n-2)R^2+(n-2)|x|^2 +|x-y|^2}{n\o_n R|x-y|^n}
&=\f{-2x(y-x)}{n\o_n R|x-y|^n} +  \f{R^2-|x|^2}{\o_n R|x-y|^n}\\
&=\f{-2xy + 2|x|^2}{n\o_n R|x-y|^n} +  n\f{R^2-|x|^2}{n\o_n R|x-y|^n}\\
&=\f{R^2-|x|^2 -2xy + 2|x|^2}{n\o_n R|x-y|^n} +  (n-1)\f{R^2-|x|^2}{n\o_n R|x-y|^n}\\
&=\f{1}{n\o_n R|x-y|^{n-2}} +  (n-1)\f{\p G(x, y)}{\p n}.
\ea\]
Therefore,
\be\ba
\div G_b(x F) &- x\cdot\g G_b (F) = \div \int_{\p\O} \f{\p G(x,y)}{\p n} y F(y)dy - x\cdot\g \int_{\p\O} \f{\p G(x,y)}{\p n} F(y) dy\\
&= \int_{\p\O} \f{y\cdot\g_x\p G(x,y)}{\p n} F(y)dy - \int_{\p\O} \f{x\cdot \p \g_x G(x,y)}{\p n} F(y) dy\\
&= \int_{\p\O} \f{(y-x)\cdot\p \g_xG(x,y)}{\p n} F(y)dy\\
&=\int_{\p\O}(\f{1}{n\o_n R|x-y|^{n-2}} +  (n-1)\f{\p G(x, y)}{\p n}) F(y) dy\\
&= (n-1)G_b(F) +\int_{\p\O} \f{F(y)}{n\o_n R|x-y|^{n-2}} dy.
\ea\ee
In view of  \eqref{ff-4}, one has
\be\ba
\f{n\ld+2\mu}{2(\mu + \ld)} F &+ \f{\mu}{\mu+\ld} g =\div G(Q) - \f{\ld}{2(\mu + \ld)}x\cdot\g G(\div f) \\
&+ \f{\ld(n-1)}{2(\mu + \ld)} G_b(F) + \f{\ld}{2(\mu + \ld)}\int_{\p\O} \f{F(y)}{n\o_n R|x-y|^{n-2}} dy.
\ea\ee
Consequently, it is easy to check that
\[\ba
\f{n\ld+2\mu}{2(\mu + \ld)} F &+ \f{\mu}{\mu+\ld} g =\div G(Q) - \f{\ld}{2(\mu + \ld)}x\cdot\g G(\div f) \\
&+ \f{\ld(n-1)}{2(\mu + \ld)} F -\f{\ld(n-1)}{2(\mu + \ld)}G(\div f) + \f{\ld}{2(\mu + \ld)}\int_{\p\O} \f{F(y)}{n\o_n R|x-y|^{n-2}} dy
\ea\]
which gives the following identity
\be\label{ff-5}
\ba
\f{2\mu+\ld}{2(\mu + \ld)} F &+ \f{\mu}{\mu+\ld} g =\div G(Q) - \f{\ld}{2(\mu + \ld)}x\cdot\g G(\div f)\\
&- \f{\ld(n-1)}{2(\mu + \ld)} G(\div f) + \f{\ld}{2(\mu + \ld)}\int_{\p\O} \f{F(y)}{n\o_n R|x-y|^{n-2}} dy.
\ea\ee

Taking $n=2$ in \eqref{ff-5} yields a special expression of F in 2D case as
\[\ba
\f{2\mu+\ld}{2(\mu + \ld)} F &+ \f{\mu}{\mu+\ld} g =\div G(Q) - \f{\ld}{2(\mu + \ld)}x\cdot\g G(\div f)\\
&- \f{\ld}{2(\mu + \ld)} G(\div f) + \f{\ld}{4\pi R(\mu + \ld)}\int_{\p\O} F(y) dy.
\ea\]
That is,
\[\ba
 \f{2\mu+\ld}{2(\mu + \ld)}F &+ \f{\mu}{\mu+\ld} g =\f{\ld}{2(\mu + \ld)}\div G(x\div f) + \f{\mu}{\mu + \ld}\div G(\g g) + \div G(f)\\
  & \qquad- \f{\ld}{2(\mu + \ld)}x\cdot\g G(\div f) - \f{\ld}{2(\mu + \ld)} G(\div f) + \f{\ld}{4\pi R(\mu + \ld)}\int_{\p\O} F(y) dy\\
&= \f{\mu}{\mu + \ld}\div G(\g g) + \div G(f) - \f{\ld}{2(\mu + \ld)} G(\div f)\\
  & \qquad- \f{\ld}{2(\mu + \ld)}\Big(\div G(x\div f)-x\cdot\g G(\div f)\Big)  + \f{\ld}{4\pi R(\mu + \ld)}\int_{\p\O} F(y) dy
\ea\]
Now one can rewrite the above identity as follows to finish the proof of Lemma \ref{lem-Farb}.
\be\ba
\f{2\mu+\ld}{2(\mu + \ld)} F & = \f{\mu}{\mu + \ld}\Big(\div G(\g g)-g\Big) + \div G(f)- \f{\ld}{2(\mu + \ld)} G(\div f)\\
  & \qquad- \f{\ld}{2(\mu + \ld)}\Big(\div G(x\div f)-x\cdot\g G(\div f)\Big)  + \f{\ld}{4\pi R(\mu + \ld)}\int_{\p\O} F(y) dy\\
    &\defeq A(g)+R(g)+B(g).
\ea\ee
where $A,R,B$ is defined as \eqref{arb}.
%\be\ba
% F & = \f{2\mu}{2\mu+\ld}\Big(\div G(\g g)-g\Big) + \f{\mu + \ld}{\mu+\ld/2}\div G(f)- \f{\ld}{2(\mu+\ld/2)} G(\div f)\\
%  & \qquad- \f{\ld}{2(\mu+\ld/2)}\Big(\div G(x\div f)-x\cdot\g G(\div f)\Big)  + \f{\ld}{4\pi R(\mu+\ld/2)}\int_{\p\O} F(y) dy.
%\ea\ee
%\be\ba
% F & = \f{2\mu}{2\mu+\ld}\Big(\div G(\g g)-g\Big) + \f{2(\mu + \ld)}{2\mu+\ld}\div G(f)- \f{\ld}{2\mu+\ld} G(\div f)\\
%  & \qquad- \f{\ld}{2\mu+\ld)}\Big(\div G(x\div f)-x\cdot\g G(\div f)\Big)  + \f{\ld}{2\pi R(2\mu+\ld)}\int_{\p\O} F(y) dy.
%\ea\ee
\endproof

Indeed, A and R enjoy the following estimates.
\begin{lemma}\label{lem-Ag}
For any $1<q<\infty$, there exists a uniform constant $\Lambda(q)$ depending only on the domain, $\lambda$ and $\mu$ such that
\be
\|A(g)+\t{g}\|_q\leq \tC(q)\|g-\t g\|_q,\quad \|\g A(g)\|_q\leq \tC(q)\|\g g\|_q.
\ee
and
\be
\|\g R(f)\|_q \leq \tC(q)\|f\|_q.
\ee
where $\bar{g}=\int_{\O}gdx$ is the average of integral over domain.
\end{lemma}
\proof
Applying  the elliptic theory to the Lame system and making use of the special structure of  A and R, the conclusion in the Lemma can be proved in a standard way.
\endproof

\section{A priori estimates of $A_1(T)$ and $A_2(T)$}

The following bounds of $\|p-\bar{p}\|_{L^q}$ will be frequently used in the estimates of $A_1(T)$, $A_2(T)$ and so on.
\begin{lemma}  \label{lem-p}
It holds that
\be
\|p-\tp(t)\|_q \leq C(\mrho,\trho) E_0^{\f 1{q\gamma}}.
\ee
for any $q\in[1, 6]$.
\end{lemma}
\proof It is easy to check that, for $1\leq q\leq 6$, it holds that
\[
\|p-p(\trho)\|_q \leq C(\mrho, \trho) E_0^{\f 1{q\gamma}}.
\]
Since
\[
|\tp(t)-p(\trho)| \leq \aint| p - p(\trho)| dx
\leq C \into| \rho^\gamma - \trho^\gamma| dx \leq C(\mrho, \trho) E_0^{\f 1\gamma},
\]
the lemma is follows.
\endproof

The next part is devoted to prove the key bootstrap argument. In order to achieve this,  we first prove the following lemma:
\begin{lemma}\label{lem-A12} It holds that
\be\label{A1}
A_1(T) \leq C E_0^{\f{1}{\gamma}} + C \int_0^T\into \sigma|\g u|^{3} d xdt+ C\int_0^T \sigma^{3}\|p-\tp\|_4^4dt.
\ee
and
\be\label{A2}
\ba
A_2(T)&\leq CE_0 + CA_1(T) + C(\mrho)A_1^{\f 32}(T)A_2^{\f 12}(T) \\
&\qquad + C(\mrho)A_2^{\f 12}(T)A_1(T)+ C\int_0^T\sigma^3\|p-\tp\|_4^4dt,
\ea\ee
provided $E_0\leq \eps_1$, where $\eps_1$ is small enough, depending on $\mu$, $\ld$, $\gamma$, $\q0$, $\trho$, $\mrho$.
\end{lemma}
\proof
{\bf Step 1.} The estimate for $A_1(T)$ is a standard procedure to derive a bound for the velocity gradient as in\cite{H1, Huang-Li-Xin}. The only subtle thing is to bound each $A_i(T)$ in terms of the $\|\rho\|_{L^q}$  rather than the $\|\rho\|_{L^\infty}$.

For any integer $m\ge 0$, taking $\sigma(t)^m\dot{u}$ as multiplier to the momentum equations with $\sigma(t)=min(1,t)$, one has
\[
\ba
{\into \sigma^{m} \rho|\dot{u}|^{2} d x} &{=\into\left(-\sigma^{m} \dot{u} \cdot \g p+\mu \sigma^{m} \Delta u \cdot \dot{u}+\ld \sigma^{m} \g \div u \cdot \dot{u}\right) d x} \\
&{\defeq M_1 +M_2 +M_3}\ea
\]

In the following, each $M_i$ will be dealt seperately as follows. First,
\[\ba M_{1} &=-\into \sigma^{m} \dot{u} \cdot \g p d x \\
&=\into\left(\sigma^{m}(\div u)_{t}(p-\tp)-\sigma^{m}(u \cdot \g u) \cdot \g p\right) d x \\
&=\f d{dt}\left(\into \sigma^{m} \div u(p-\tp) d x\right) -m \sigma^{m-1} \sigma^{\prime} \into \div u(p-\tp) d x \\
    &\quad+\into \sigma^{m}\left(p' \rho(\div u)^{2}-p(\div u)^{2}+ p \p_{i} u^{j} \p_{j} u^{i}\right) d x\\
&\leq \f{d}{dt}\left(\into \sigma^{m} \div u(p-\tp) d x\right) + C\sigma'\into |\g u||p- \tp| dx\\
    &\qquad + C\into\sigma^m|\g u|^2|p- \tp|dx + C\E_0\|\g u\|_2^2,
%&\leq \left(\into \sigma^{m} \div u(p-\tp) d x\right)_{t} + C\sigma'\|\g u\|_2\|\max(\rho,\trho)^{\gamma-1}\|_4\|\rho- \trho\|_4  \\
%  &\qquad+ C\sigma^m\big(\into|\g u|^3dx)\big)^{\f 23}\big(\into|p- \tp|^3dx\big)^{\f 13}\\
%&\leq \left(\into \sigma^{m} \div u(p-\tp) d x\right)_{t} + C(\trho)\sigma'\|\g u\|_2\|\rho- \trho\|_\gamma^{\f 14}\|\rho- \trho\|_{\f{4\gamma}{\gamma-1}}^{\f 34}  \\
%  &\qquad+ C\sigma^m\big(\into|\g u|^3dx)\big)^{\f 23}\big(\into|p- \tp|^3dx\big)^{\f 13}\\
%&\leq\left(\into \sigma^{m} \div u(p-\tp) d x\right)_{t} + C(\trho)\sigma'\|\g u\|_2\|\rho- \trho\|_\gamma^{\f 14}
%    + C(\trho)E_0^{\f 1{3\gamma}}\sigma^m\big(\into|\g u|^3dx)\big)^{\f 23}\\
%&\leq \left(\into \sigma^{m} \div u(p-\tp) d x\right)_{t} + C(\trho)\|\g u\|_2^2 + C(\trho)\sigma'E_0^{\f 1{2\gamma}} + C(\trho)E_0^{\f 1{3\gamma}}\sigma^m\big(\into|\g u|^3dx)\big)^{\f 23}
\ea\]
and
\[
\begin{aligned}
M_{2} &=\into \mu \sigma^{m} \Delta u \cdot \dot{u} d x \\
&=-\frac{\mu}{2}\f{d}{dt}\left(\sigma^{m}\|\g u\|_{2}^{2}\right)+\frac{\mu m}{2} \sigma^{m-1} \sigma^{\prime}\|\g u\|_{2}^{2}
    -\mu \sigma^{m} \into \p_{i} u^{j} \p_{i}\left(u^{k} \p_{k} u^{j}\right) d x \\
& \leq-\frac{\mu}{2}\f{d}{dt}\left(\sigma^{m}\|\g u\|_{2}^{2}\right)+C m \sigma^{m-1}\|\g u\|_{2}^{2}+C \into \sigma^{m}|\g u|^{3} d x
\end{aligned}
\]
$M_3$ can be treated in a similar way as
\[
\begin{aligned}
M_{3}=&-\frac{\ld}{2}\f{d}{dt}\left(\sigma^{m}\|\div u\|_{2}^{2}\right)+\frac{m\ld}{2} \sigma^{m-1}\|\div u\|_{2}^{2} \\
 &-\ld \sigma^{m} \into \div u \div(u \cdot \g u) d x \\
  \leq &-\frac{\ld}{2}\f{d}{dt}\left(\sigma^{m}\|\div u\|_{2}^{2}\right)
  +C m \sigma^{m-1}\|\g u\|_{2}^{2}+C \into \sigma^{m}|\g u|^{3} d x
   \end{aligned}
\]
Combing all the above estimates together yields that
\be\label{uu-1}
\ba
\f{d}{dt}\left(\sigma^{m} B(t)\right)&+\into \sigma^{m} \rho|\dot{u}|^{2} d x \leq C m \sigma^{m-1}\|\g u\|_{2}^{2} +C \into \sigma^{m}|\g u|^{3} d x\\
&+ C\sigma'\into |\g u||p- \tp| dx + C\into\sigma^m|\g u|^2|p- \tp|dx,
\ea
\ee
with
\[
\ba
B(t) & \defeq \frac{\mu}{2}\|\g u\|_{2}^{2}+\frac{\ld}{2}\|\div u\|_{2}^{2}-\into \div u(p-\tp) d x \\
 & \geq \frac{\mu}{2}\|\g u\|_{2}^{2}+\frac{\ld}{2}\|\div u\|_{2}^{2}-C(\mrho) E_{0}^{\f 1 {2\gamma}}\|\div u\|_{2} \\
 & \geq \frac{\mu}{4}\|\g u\|_{2}^{2}+\frac{\ld}{4}\|\div u\|_{2}^{2}-C(\mrho) E_{0}^{\f 1\gamma}.
\ea
\]
where  Lemma \ref{lem-p} has been used.

Taking $m=1$ in \eqref{uu-1} gives
\[\ba
\int_0^T & \sigma'\into |\g u||p- \tp| dxdt\leq C\int_0^{\sigma(T)}\|\g u\|_2\|p- \tp\|_2dt\\
&\leq C(\mrho)E_0^{\f 1{2\gamma}}\int_0^{\sigma(T)}\|\g u\|_2dt
\leq C(\mrho)E_0.
\ea\]
It remains to deal with the last term in \eqref{uu-1}. We decompose the integral into one over short time $(0,\sigma(T))$ and the other over long time $(\sigma(T),T)$.

 Using Lemma \ref{lem-du} with $q=6$ and Lemma \ref{lem-p}, the short time integral is bounded as follows:
\[\ba
\int_0^{\sigma(T)}\into \sigma|\g u|^2|p-\tp|dxdt &\leq \int_0^{\sigma(T)}\sigma\|\g u\|_2\|\g u\|_6\|p-\tp\|_3dt\\
%&\leq \int_0^{\sigma(T)}\sigma\|\g u\|_3^2\|p-\tp\|_\f{3}{2}dt \leq C(\mrho)E_0^{\f{2}{3\gamma}}\int_0^{\sigma(T)}\sigma\|\g u\|_2\|\g u\|_6dt\\
%&\leq C(\mrho)E_0^{\f{2}{3\gamma}}\int_0^{\sigma(T)}\sigma\|\g u\|_2(\|\g w\|_6 +\|\g v\|_6)dt\\
%&\leq C(\mrho)E_0^{\f{2}{3\gamma}}\int_0^{\sigma(T)}\sigma\|\g u\|_2(\|\g w\|_6 +\| p-\tp\|_6)dt\\
%&\leq C(\mrho)E_0^{\f{2}{3\gamma}}\int_0^{\sigma(T)}\sigma\|\g u\|_2(\|\g^2 w\|_{\f 32} +\| p-\tp\|_6)dt\\
%&\leq C(\mrho)E_0^{\f{2}{3\gamma}}\int_0^{\sigma(T)}\sigma\|\g u\|_2(\|\rho \dot u\|_{\f 32} +\| p-\tp\|_6)dt\\
%&\leq C(\mrho)E_0^{\f{2}{3\gamma}}\int_0^{\sigma(T)}\sigma\|\g u\|_2\|\rho^{\f12} \dot u\|_2\|\rho^{\f 12}\|_6  dt + C(\mrho)E_0^{\f{2}{3\gamma}}\\
&\leq C(\mrho)E_0^{\f{1}{3\gamma}}\int_0^{\sigma(T)}\sigma\|\g u\|_2 \big(\|\rho^{\f12} \dot u\|_2\|\rho^{\f 12}\|_6 +\| p-\tp\|_6\big) dt\\
&\leq C(\mrho)E_0^{\f{1}{3\gamma}}\int_0^{\sigma(T)}\sigma(\|\g u\|_2^2 + \|\rho^{\f12} \dot u\|^2_2) dt + C(\mrho)E_0\\
&\leq C(\mrho)E_0 + C_1(\mrho)E_0^{\f{1}{3\gamma}}A_1(T).
\ea\]
The integral over long time integral can be aslo bounded  as
\[\ba
\int_\sigma^T\into &\sigma|\g u|^2|p-\tp|dxdt \leq \int_\sigma^T \|\g u\|_2\|\g u\|_4\|p-\tp\|_4dt\\
&\leq \int_\sigma^T \|\g u\|_2(\|\g w\|_4 +\|p-\tp\|_4)\|p-\tp\|_4dt\\
&\leq \int_\sigma^T \|\g u\|_2\|\g^2 w\|_{\f 43}\|p-\tp\|_4dt + \int_\sigma^T \|\g u\|_2\|p-\tp\|_4^2dt\\
&\leq \int_\sigma^T \|\g u\|_2\|\rho \dot u\|_{\f 43}\|p-\tp\|_4dt + \int_\sigma^T \|\g u\|_2^2 + \|p-\tp\|_4^4dt\\
&\leq \int_\sigma^T \|\g u\|_2\|\rho^{\f 12} \dot u\|_2\|\rho^{\f 12}\|_4\|p-\tp\|_4dt + \int_\sigma^T \|p-\tp\|_4^4dt + CE_0\\
&\leq C_2(\mrho)E_0^{\f 1{4\gamma}}\int_\sigma^T \|\g u\|_2^2 + \|\rho^{\f 12} \dot u\|^2_2dt + \int_\sigma^T \|p-\tp\|_4^4dt + CE_0\\
&\leq C(\mrho)E_0 + C_2(\mrho)E_0^{\f 1{4\gamma}}A_1(T) + \int_\sigma^T \|p-\tp\|_4^4dt.
\ea\]
Therefore, integrating \eqref{uu-1} over $(0,T)$ yields
\be\label{uu-2}
\ba
\sigma B(T)&+\int_0^T\into \sigma \rho|\dot{u}|^{2} d x \leq C(\mrho) E_{0}^{\f 1\gamma} +C \int_0^T\into \sigma|\g u|^{3} d xdt\\
&\qquad + C_1(\mrho)E_0^{\f{1}{3\gamma}}A_1(T) + \int_\sigma^T \|p-\tp\|_4^4dt + CE_0.
\ea
\ee
Set $\eps_1=\big(\f{\mu}{8C_1(\mrho)}\big)^{3\gamma}$ in \eqref{uu-2}. Then the bound for $A_1(T)$ in
\eqref{A1} follows immediately provided that $E_0\leq \eps_1$.

{\bf Step 2.} Here we will use an idea due to Hoff to deal with $A_2(T)$, i.e, for any integer $m\ge 0$, operating $\sigma^m\dot{u}^j[\p/\p t+div(u\cdot)]$ on the $j_{th}$ part of the momentum equation, summing with respect to $j$, and integrating the resulting equation over $\O$, one obtains after integration by parts
\[\ba
\f{d}{dt}\left(\frac{\sigma^{m}}{2} \into \rho|\dot{u}|^{2} d x\right)&-\frac{m}{2} \sigma^{m-1} \sigma^{\prime} \into \rho|\dot{u}|^{2} d x\\
&=-\into \sigma^{m} \dot{u}^{j}\left[\p_{j} p_{t}+\div\left(\p_{j} p u\right)\right] d x\\
&+\mu \into \sigma^{m} \dot{u}^{j}\left[\triangle u_{t}^{j}+\div\left(u \Delta u^{j}\right)\right] d x\\
&+\ld \into \sigma^{m} \dot{u}^{j}\left[\p_{t} \p_{j} \div u+\div\left(u \p_{j} \div u\right)\right] d x\\
&\defeq N_1 +N_2 + N_3.
\ea\]
\def\ol#1{\overline{#1}}
\def\tp{\overline{p}}
It follows from the mass equation that
\[\begin{aligned}
 N_{1} &=-\into \sigma^{m} \dot{u}^{j}\left[\p_{j} p_{t}+\div\left(\p_{j} p u\right)\right] d x \\
 &=\into \sigma^{m}\Big[-\gamma(p-\tp)  \div u \p_{j} \dot{u}^{j}+\p_{k}\left(\p_{j} \dot{u}^{j} u^{k}\right) (p-\tp) \\
 &\qquad - (p-\tp)\p_{j}\left(\p_{k} \dot{u}^{j}  u^{k}\right) -\gamma\tp\p_j\dot u_j\div u \Big] d x \\
&\leq C\into\sigma^m |p -\tp||\g u||\g\dot u|dx + C\E_0\into\sigma^m |\g u| |\g\dot u|dx
%  & \leq C\sigma^{m}\|\g \dot{u}\|_{2}(\|\g u\|_{L^{4}}\|p-\tp\|_4 + \|\g u\|_{L^{4}}\|\rho p' -\ol{\rho p'}\|_4 + \|\g u\|_2) \\
%   & \leq \delta \sigma^{m}\|\g \dot{u}\|_{2}^{2} + C(\delta)\sigma^{m}\|\g u\|_2^2 +C(\overline{\rho}, \delta) E_0\sigma^{m}\|\g u\|_{L^{4}}^{2}
\end{aligned}\]
One can integrate by part to obtain the following bound for $N_2$.
\[\begin{aligned}
N_{2} &=\mu \into \sigma^{m} \dot{u}^{j}\left[\triangle u_{t}^{j}+\div\left(u \Delta u^{j}\right)\right] d x \\
 &=-\mu \into \sigma^{m}\left[|\g \dot{u}|^{2}+\p_{i} \dot{u}^{j} \p_{k} u^{k} \p_{i} u^{j}
-\p_{i} \dot{u}^{j} \p_{i} u^{k} \p_{k} u^{j}-\p_{i} u^{j} \p_{i} u^{k} \p_{k} \dot{u}^{j} \right] d x\\
 & \leq-\frac{3 \mu}{4} \into \sigma^{m}|\g \dot{u}|^{2} d x+C \into \sigma^{m}|\g u|^{4} d x
 \end{aligned}\]
Finally, $N_3$ can be handled in a same way as for $N_2$.
\[
N_{3} \leq-\frac{\lambda}{2} \into \sigma^{m}(\div \dot{u})^{2} d x+C \into \sigma^{m}|\g u|^{4} d x
\]
Combing all the above together yields
\be\label{mdotu}
\ba
\f{d}{dt}\Big(\sigma^{m}& \into \rho|\dot{u}|^{2} d x\Big)+\mu \into \sigma^{m}|\g \dot{u}|^{2} d x+ \lambda\into \sigma^{m}(\div \dot{u})^{2} d x\\
&\leq m \sigma^{m-1} \sigma^{\prime} \into \rho|\dot{u}|^{2} d x + C\into\sigma^m |p -\tp||\g u||\g\dot u|dx \\
    &\qquad + C\into\sigma^m |\g u| |\g\dot u|dx  + C \into \sigma^{m}|\g u|^{4} d x\\
&\leq m \sigma^{m-1} \sigma^{\prime} \into \rho|\dot{u}|^{2} d x+ C\into\sigma^m |p -\tp|^2|\g u|^2 dx \\
    &\qquad+ C\into\sigma^m |\g u|^2dx  + C \into \sigma^{m}|\g u|^{4} d x
\ea\ee
 To finish the estimate of $A_2(T)$, one takes $m=3$ in \eqref{mdotu} to get
\[\ba
\sigma^3 \into \rho|\dot{u}|^{2} d x& + \mu \int_0^T\into \sigma^3|\g \dot{u}|^{2} d xdt+ \lambda\int_0^T\into \sigma^3(\div \dot{u})^{2} d xdt\\
&\leq C\int_0^{\sigma(T)} \into\sigma^2 \rho|\dot{u}|^{2} d x + C\int_0^T\into\sigma^3 |p -\tp|^2|\g u|^2dx \\
    &\qquad + C\int_0^T\into\sigma^3 |\g u|^2dx  + C \int_0^T\into \sigma^{3}|\g u|^{4} d x\\
&\leq C\int_0^{\sigma(T)} \into\sigma^2 \rho|\dot{u}|^{2} d x + C\int_0^T\into\sigma^3 |p -\tp|^4dx + CE_0  + C \int_0^T\into \sigma^{3}|\g u|^{4} d x
\ea\]

It remains to treat $\|\nabla u\|_{4}$. Direct computation lead to
\[\ba
\|\g u\|_4^4 &\leq C\|\g w\|_4 + C\|\g v\|_4\\
&\leq C\|\g w\|_2\|\g w\|_6^3 + C\|p-\tp\|_4^4\\
&\leq C(\|\g u\|_2 + \|p-\tp\|_2)\|\g^2 w\|_{\f 32}^3 + C\|p-\tp\|_4^4\\
&\leq C(\|\g u\|_2 + \|p-\tp\|_2)\|\rho \dot u\|_{\f 32}^3 + C\|p-\tp\|_4^4\\
&\leq C(\|\g u\|_2 + \|p-\tp\|_2)\|\rho^{\f 12} \dot u\|_{2}^3\|\rho^{\f 12}\|_6 + C\|p-\tp\|_4^4\\
&\leq C(\mrho)(\|\g u\|_2 + \|p-\tp\|_2)\|\rho^{\f 12} \dot u\|_{2}^3 + C\|p-\tp\|_4^4\\
&\leq C(\mrho)\|\g u\|_2\|\rho^{\f 12}\dot u\|_{2}^3 + C(\mrho)\|\rho^{\f 12}\dot u\|_{2}^3  + C\|p-\tp\|_4^4.
\ea\]
Hence,
\[\ba
\int_0^T\sigma^3\|\g u\|_4^4dt &\leq C(\mrho)\int_0^T\sigma^3\|\g u\|_2\|\rho^{\f 12} \dot u\|_{2}^3dt + C(\mrho)\int_0^T\sigma^3\|\rho^{\f 12} \dot u\|_{2}^3 dt \\
    &\qquad + C\int_0^T\sigma^3\|p-\tp\|_4^4dt\\
&\leq C(\mrho)\sup\limits_{0\leq t\leq T}\sigma^{\f12}\|\g u\|_2 \sup\limits_{0\leq t\leq T}\sigma^{\f 32}\|\rho^{\f 12} \dot u\|_{2}\int_0^T\sigma\|\rho^{\f 12} \dot u\|_{2}^2dt \\
    &\qquad+ C(\mrho)\sup\limits_{0\leq t\leq T}\sigma^{\f 32}\|\rho^{\f 12} \dot u\|_{2}\int_0^T\sigma^{\f 32}\|\rho^{\f 12} \dot u\|_{2}^2 dt
        + C\int_0^T\sigma^3\|p-\tp\|_4^4dt\\
&\leq C(\mrho)A_1^{\f 32}(T)A_2^{\f 12}(T) + C(\mrho)A_2^{\f 12}(T)A_1(T)+ C\int_0^T\sigma^3\|p-\tp\|_4^4dt.
\ea\]
Clearly, the first term in \eqref{mdotu} is bounded by $A_1(T)$
\[
\int_0^{\sigma(T)}\sigma^2\into\rho|\dot u|^2dx dt\leq CA_1(T).
\]
This finishes the proof of \eqref{A2}.
\endproof

\section{ Short time a priori estimates of $\g u$}
Next, the following lemma will play important roles in the estimates on the uniform upper bound of the density for small time.
\begin{lemma}\label{lem-tgu}
It holds that, for $\beta\in(\f 12, 1]$,
\be\label{tgu1}
\sup _{0 \leq t \leq \sigma(T)} \sigma^{1-\beta}\|\g u\|_{2}^{2}+\int_{0}^{\sigma(T)} \sigma^{1-\beta} \into \rho|\dot{u}|^{2} d x d t \leq C(\mrho, M),
\ee
and
\be\label{tgu2}
\sup _{0 \leq t \leq \sigma(T)} \sigma^{2-\beta} \into \rho|\dot{u}|^{2} d x+\int_{0}^{\sigma(T)} \sigma^{2-\beta} \into|\g \dot{u}|^{2} d x d t \leq C(\mrho, M),
\ee
provided that $E_0<\eps_2$,
with $\eps_2$ is small enough depending on $\mu$, $\ld$, $\gamma$, $q_0$, $\trho$, $\mrho$, $\beta$, $\d_0$, $M$.
\end{lemma}
\proof As in \cite{Hof3}, define $\omega_1$ and $\omega_2$ to be the solutions to the following problems
\be\label{w1-1}
\left\{\ba
&\rho \dot w_1 = \mu\lap w_1 +\ld\g\div w_1,\\
&w_1(x,0) = w_{10}
\ea\right.\ee
and
\be\label{w2-1}
\left\{\ba
&\rho \dot w_2 + \g p = \mu\lap w_2 +\ld\g\div w_2,\\
&w_2(x,0) = w_{20}
\ea\right.\ee
respectively.

{\bf Step 1.} Standard energy estimates show that
\be
\sup _{0 \leq t \leq \sigma(T)} \into \rho\left|w_{1}\right|^{2} d x+\int_{0}^{\sigma(T)} \into\left|\g w_{1}\right|^{2} d x d t \leq \|\rho_0\|_{\infty}\into\left|w_{10}\right|^{2} d x
\ee
and
\be
\sup _{0 \leq t \leq \sigma(T)} \into \rho|w_{2}|^{2} d x+\int_{0}^{\sigma(T)} \into|\g w_{2}|^{2} d x d t \leq \into\rho_0|w_{20}|^{2} d x + C(\mrho) E_{0}^{\f 1{\gamma}}.
\ee
{\bf Step 2.} Estimate of $w_1$. Multiplying \eqref{w1-1} by $\omega_{1t}$ and integrating the resulting equality over $\O$, one gets
\[\ba
\frac{1}{2}\f{d}{dt}\Big(\mu\|\g w_{1}\|_{2}^{2}&+\ld\|\div w_{1}\|_{2}^{2}\Big) + \into \rho|\dot{w}_{1}|^{2} d x
    =\into (\rho u \cdot \g w_{1})\dot{w_{1}}d x\\
&\leq C \|\rho^{\f 12}\dot w_1\|_2 \|\rho^{\f 13}u\|_3 \|\g w_1\|_{12}\|\rho^{\f 12}\|_{12}
\leq C_1(\mrho)E_0^{\f{\d_0}{3}} \|\rho^{\f 12}\dot w_1\|_2^2.
\ea\]
Therefore
\[\ba
&\sup _{0 \leq t \leq \sigma(T)}\|\g w_{1}\|_{2}^{2}+\int_{0}^{\sigma(T)} \into \rho|\dot{w}_{1}|^{2} d x d t \leq C(\mrho)\|\g w_{10}\|_{2}^{2},\\
&\sup _{0 \leq t \leq \sigma(T)} \sigma\|\g w_{1}\|_{2}^{2}+\int_{0}^{\sigma(T)} \sigma \into \rho|\dot{w}_{1}|^{2} d x d t \leq C(\mrho, M)\|w_{10}\|_{2}^{2}.
\ea\]
provided that $\eps_2\leq \Big(\f 1{2C_1(\mrho)}\Big)^{\f{3}{\d_0}}$.

Then, by the standard interpolation theory, for any $\th\in [0,1]$, one has
\be
\sup _{0 \leq t \leq \sigma(T)} \sigma^{1-\theta}\|\g w_{1}\|_{2}^{2}+\int_{0}^{\sigma(T)} \sigma^{1-\theta} \into \rho|w_{1}|^{2} d x d t \leq C(\mrho)\|w_{10}\|_{\dot{H}^{\theta}}^{2}.
\ee

{\bf Step 3.} Estimate of $w_2$. Recall that $v$ solves
\[\left\{\ba
&\mu\lap v + \ld\g\div v = \g p,\\
& v\big|_{\p\O} = 0.
\ea\right.,
\]
Then
\[
\sup\limits_{0 \leq t \leq \sigma(T)}\|\g v\|_2\leq C\sup\limits_{0 \leq t \leq \sigma(T)}\|p-\tp\|_2 \leq C(\mrho)E_0^{\f 1{2\gamma}}.
\]
Set $w_3 = w_2-v$. Then
\[\left\{\ba
&\mu\lap w_3 + \ld\g\div w_3 = \rho\dot w_2\\
& w_3\big|_{\p\O} = 0.
\ea\right.
\]
and
\[
\|\g w_{30}\|_2 \leq \|\g w_{20}\|_2 + \|\g v_{0}\|\leq \|\g w_{20}\|_2 +\|p(\rho_0)-\tp(0)\|_2 \leq \|\g w_{20}\|_2 + C(\mrho)E_0^{\f 1{2\gamma}}.
\]

Then
\[\ba
\frac{1}{2}\f{d}{dt}\Big(\mu\|\g w_{3}\|_{2}^{2}&+\ld\|\div w_{3}\|_{2}^{2}\Big) + \into \rho|\dot{w}_{2}|^{2}  d x\\
    &= \into (\rho u \cdot \g w_{3})\dot{w_{2}} d x + \into \rho \dot v\cdot\dot{w_{2}} d x \\
&\leq \|\rho^{\f 12}\dot w_2\|_2 \|\rho^{\f 13}u\|_3 \|\g w_3\|_6 + \|\rho^{\f 12}\dot w_2\|_2\|\rho^{\f 12}\dot v\|_2\\
&\leq \|\rho^{\f 12}\dot w_2\|_2 \|\rho^{\f 13}u\|_3 \|\g^2 w_3\|_{\f 3 2} + \|\rho^{\f 12}\dot w_2\|_2\|\rho^{\f 12}\dot v\|_2\\
&\leq C_2(\mrho)E_0^{\f{\d_0}{3}} \|\rho^{\f 12}\dot w_2\|_2^2+ \f 14\|\rho^{\f 12}\dot w_2\|_2^2 + \|\rho^{\f 12}\dot v\|_2^2.
\ea\]
If $C_2(\mrho)E_0^{\f{\d_0}{3}}\leq \f 12$, then one has
\[\ba
\f{d}{dt}\Big(\mu\|\g w_{3}\|_{2}^{2}&+\ld\|\div w_{3}\|_{2}^{2}\Big) + \into \rho|\dot{w}_{2}|^{2}  d x
\leq C\|\rho^{\f 12}\dot v\|_2^2.
\ea\]
Noting that
\[
\|\rho^{\f 12}\dot v\|_2 \leq \|\rho^{\f 12}\p_t v\|_2 + \|\rho^{\f 12}u\cdot\g v\|_2
\]
and
\[\ba
\|\rho^{\f 12}\p_t v\|_2 &\leq C(\mrho) \|L^{-1}\g p_t\|_3\leq C(\mrho) \|L^{-1}\g (\div (pu) + (\gamma-1)p\div u)\|_3\\
&\leq C(\mrho) \big(\|pu\|_3 + \|\g L^{-1}\g(p\div u)\|_{\f6 5}\big) \leq C(\mrho) \big(\|\g u\|_2 + \|p\div u\|_{\f6 5}\big)\\
&\leq C(\mrho) \|\g u\|_2,
\ea\]
and
\[
\|\rho^{\f 12}u\cdot\g v\|_2 \leq \|\rho^{\f 13}u\|_3 \|\rho^{\f 16}\|_{12}\|\g v\|_{12}\leq C(\mrho) E_0^{\f {\d_0}{3}}.
\]
Then
\[\ba
\f{d}{dt}\Big(\mu\|\g w_{3}\|_{2}^{2}&+\ld\|\div w_{3}\|_{2}^{2}\Big) + \into \rho|\dot{w}_{2}|^{2}  d x
\leq C(\mrho) \|\g u\|_2^2 + C(\mrho) E_0^{\f {2\d_0}3}.
\ea\]
Hence,
\[\ba
\sup\limits_{0\leq t\leq\sigma(T)}&\|\g w_3\|_2^2 + \int_0^{\sigma(T)}\into\rho |\dot w_2|^2dxdt \leq \|\g w_{30}\|_2^2 + C(\mrho) E_0 + C(\mrho) E_0^{\f {2\d_0}3}\\
& \leq \|\g w_{20}\|_2^2 + C(\mrho)E_0^{\f 1{2\gamma}} + C(\mrho) E_0^{\f {2\d_0}3}.
\ea\]
So
\[\ba
\sup\limits_{0\leq t\leq\sigma(T)}&\|\g w_2\|_2^2 + \int_0^{\sigma(T)}\into\rho |\dot w_2|^2dxdt
    \leq \|\g w_{20}\|_2^2 + C(\mrho)E_0^{\f 1{2\gamma}} + C(\mrho) E_0^{\f {2\d_0}3}.
\ea\]
Taking $w_{10} = u_0(x)$, $w_{20} = 0$, one then gets that $u=w_1 + w_2$. Therefore, \eqref{tgu1} is proved, provided that
$E_0\leq \eps_2=\min\Big\{\Big(\f 1{2C_1(\mrho)}\Big)^{\f{3}{\d_0}},\ \Big(\f 1{2C_2(\mrho)}\Big)^{\f{3}{\d_0}}\Big\}$.

{\bf Step 4.}  To prove \eqref{tgu2}, taking $m=2-\beta$ in \eqref{mdotu}, one gets
\[\ba
\sup _{0 \leq t \leq \sigma(T)} &\sigma^{2-\beta} \into \rho|\dot{u}|^{2} d x+\int_{0}^{\sigma(T)} \sigma^{2-\beta} \into|\g \dot{u}|^{2} d x d t
    \leq \int_0^{\sigma(T)} \sigma^{1-\beta} \into \rho|\dot{u}|^{2} d xdt \\
    &\quad + C(\mrho)\int_0^{\sigma(T)} \sigma^{2-\beta}\into |\g u|^4dxdt + CE_0\\
&\leq C(\mrho, M) + C(\mrho)\int_0^{\sigma(T)} \sigma^{2-\beta}\into |\g u|^4dxdt.
\ea\]
Note that
\[\ba
\int_0^{\sigma(T)} &\sigma^{2-\beta}\into |\g u|^4dxdt \leq \int_0^{\sigma(T)} \sigma^{2-\beta}\|\g u\|_2\|\g u\|_6^3dt\\
    &\leq C\int_0^{\sigma(T)} \sigma^{2-\beta}\|\g u\|_2(\|\rho^{\f 12}\|_6^3\|\rho^{\f 12}\dot u\|_2^3 + \|p-\tp\|_6^3)dt\\
    &\leq C(\mrho)\int_0^{\sigma(T)} \sigma^{2-\beta}\|\g u\|_2 \|\rho^{\f 12}\dot u\|_2^3dt + C(\mrho)\\
    &\leq C(\mrho)\int_0^{\sigma(T)} \sigma^{\beta-\f 12}(\sigma^{1-\beta}\|\g u\|_2^2)^{\f 12} (\sigma^{1-\beta}\|\rho^{\f 12}\dot u\|_2^2) (\sigma^{2-\beta}\|\rho^{\f 12}\dot u\|_2^2)^{\f 12}dt + C(\mrho)\\
    &\leq C(\mrho, M) \sup\limits_{0\leq t\leq\sigma(T)}\big(\sigma^{2-\beta}\|\rho^{\f 12}\dot u\|_2^2\big)^{\f 12}+ C(\mrho).
\ea\]
This proves \eqref{tgu2}.
\endproof

Next, we will bound the auxillary function $A_3(T)$.  It can be done in a similary way in  \cite{Huang-Li-Xin}.
\begin{lemma}\label{lem-ru3}
It holds that
\be
A_3(\sigma(T))=\sup\limits_{0\leq t\leq\sigma(T)}\into \rho |u|^3dx \leq E_0^{\d_0},
\ee
provided that $E_0\leq \eps_3$,
with $\eps_3$ is small enough depending on $\mu$, $\ld$, $\gamma$, $q_0$, $\trho$, $\mrho$, $\beta$, $\d_0$, $M$.
\end{lemma}
\proof In a similar way in \cite{Huang-Li-Xin}, one can get
\[\ba
\f{d}{dt}\into\rho |u|^3 dx &\leq C\into|u||\g u|^2 + C\into|p-\tp||u||\g u| dx\\
&\leq C\|u\|_6\|\g u\|_3\|\g u\|_2 + C \|p-\tp\|_6\|u\|_3\|\g u\|_2\\
&\leq C\|\g u\|_2^2\|\g u\|_3 + C(\mrho) \|\g u\|_2^2
    \leq C\|\g u\|_2^2\|\g u\|_2^{\f 12}\|\g u\|_6^{\f 12} + C(\mrho) \|\g u\|_2^2\\
&\leq C(\mrho)\|\g u\|_2^{\f 52}(\|\rho^{\f 12}\dot u\|_2^{\f 12} + \|p-\tp\|_6^{\f 12}) + C(\mrho) \|\g u\|_2^2\\
&\leq C(\mrho)\|\g u\|_2^{\f 52}\|\rho^{\f 12}\dot u\|_2^{\f 12} + C(\mrho)\|\g u\|_2^{\f 52} + C(\mrho) \|\g u\|_2^2.
\ea\]
Since $\beta\in (\f 12, 1]$,
\[
\frac{2\left(3-4 \delta_{0}\right)(1-\beta)}{3-8 \delta_{0}}=1-\f{\beta(2\beta-1)}{8-7\beta}<1,
\]
So
\[\ba
\int_0^{\sigma(T)}&\|\g u\|_2^{\f 52}\|\rho^{\f 12}\dot u\|_2^{\f 12}dt
    =\int_0^{\sigma(T)} \sigma^{-(1-\beta)[\f 32-2\d_0]}\big(\sigma^{1-\beta}\|\g u\|_2^2\big)^{\f 54-2\d_0} \big(\sigma^{1-\beta}\|\rho^{\f 12}\dot u\|_2^2\big)^{\f 14} \|\g u\|_2^{4\d_0} dt \\
&\leq C(\mrho, M)\int_0^{\sigma(T)} \sigma^{-(1-\beta)[\f 32-2\d_0]}\big(\sigma^{1-\beta}\|\rho^{\f 12}\dot u\|_2^2\big)^{\f 14} \|\g u\|_2^{4\d_0} dt\\
&\leq C(\mrho, M)\Big(\int_{0}^{\sigma(T)} \sigma^{-\frac{2(3-4 \delta_{0})(1-\beta)}{3-8 \delta_{0}}} d t\Big)^{(3-8 \delta_{0}) / 4}
    \Big(\int_{0}^{\sigma(T)}\|\g u\|_{2}^{2} d t\Big)^{2 \delta_{0}}\\
&\leq C(\mrho) E_0^{2\d_0},
\ea\]
Similarly,
\[\ba
\int_0^{\sigma(T)} &\|\g u\|_2^{\f 52} dt \leq \int_0^{\sigma(T)}\sigma^{-3(1-\beta) / 4}\Big(\sigma^{1-\beta}\|\g u\|_{2}^{2}\Big)^{3 / 4}\|\g u\|_{2}dt\\
&\leq C(\mrho, M) \left(\int_{0}^{\sigma(T)} \sigma^{-3(1-\beta) / 2} d t\right)^{1 / 2}\left(\int_{0}^{\sigma(T)}\|\g u\|_{2}^{2} d t\right)^{1 / 2}\\
&\leq C(\mrho, M) E_0^{\f 12}.
\ea\]
Note that
\[
\into \rho_{0}\left|u_{0}\right|^{3} d x \leq C(\mrho)\left(\into \rho_{0}\left|u_{0}\right|^{2} d x\right)^{2\d_0}\left\|u_{0}\right\|_{\dot{H}^{\beta}}^{3-4\d_0}
\leq C(\mrho, M)E_0^{2\d_0}.
\]
Therefore,
\[
\sup\limits_{0\leq t\leq\sigma(T)}\into\rho |u|^3dx \leq C_1(\mrho, M)E_0^{2\d_0}\leq E_0^{\d_0},
\]
provided that $E_0\leq \eps_3$ with $\eps_3=\min\Big\{\eps_1,\ \eps_2,\ \Big(\f 1{2C_1(\mrho, M)}\Big)^{\f {1}{\d_0}}\Big\}$.
\endproof

\section{Long time a priori estimates of $\g u$}
We now turn to the long time a priori estimates. We start with the following lemma.

\begin{lemma}\label{lem-A123}
It holds that
\be\label{A12est}
A_1(T)+A_2(T)\leq E_0^{\f 1{2\gamma}}.
\ee
and
\be\label{A3est}
A_3(T) \leq C(\mrho)\big(E_0^{\d_0} + E_0^{\f 1{2\gamma}}\big).
\ee
provided that $\f{2\mu}{2\mu+\ld}\leq \f 1{4 \tC^{\f 14}(4)}$ and $E_0\leq\eps_4$
with $\eps_4$ is small enough depending on $\mu$, $\ld$, $\gamma$, $q_0$, $\trho$, $\mrho$, $\beta$, $\d_0$, $M$.
\end{lemma}
\proof Clearly, \eqref{A3est} can be easily deduced from (6.1) and lemma \ref{lem-ru3}. So it  suffices to prove \eqref{A12est}.

Using \eqref{H1} and lemma \ref{lem-A12}, one has
\[\ba
A_1(T) + A_2(T)\leq C E_0^{\f{1}{\gamma}} + C \int_0^T\into \sigma|\g u|^{3} d xdt+ C\int_0^T \sigma^{3}\|p-\tp\|_4^4dt.
\ea\]
provided that $E_0\leq\eps_1$. So we need to estimate the last two terms in right hand side in above inequality.

{\bf Step 1.} Estimate of $\int_0^T\into \sigma|\g u|^{3} d xdt$.
Since, for $t \in [\sigma(T),\ T]$, one has
\[\ba
\|\g u\|_4^4 &\leq C\|\g w\|_4^4 + \|p-\tp\|_4^4 \leq C\|\g^2 w\|_{\f {4}{3}}^4 + \|p-\tp\|_4^4\\
&\leq \|\rho^{\f 12}\dot u\|_2^{4}\|\rho^{\f 12}\|_4^4 + \|p-\tp\|_4^4\\
&\leq C(\mrho )A_2(T)\|\rho^{\f 12}\dot u\|_2^{2}+ \|p-\tp\|_4^4.
\ea\]
Then
\[\ba
\int_{\sigma(T)}^T\into &\sigma|\g u|^{3} d xdt \leq C\int_{\sigma(T)}^T \|\g u\|_2^2 + \|\g u\|_4^4dt \leq C E_0 + \int_{\sigma(T)}^T\sigma^3\|\g u\|_4^4dt\\
&\leq C E_0 + C(\mrho )A_2(T)\int_{\sigma(T)}^T\sigma^3\|\rho^{\f 12}\dot u\|_2^{2}dt + \int_{\sigma(T)}^T\sigma^3\|p-\tp\|_4^4dt\\
&\leq C E_0 + C(\mrho )A_2(T)A_1(T) + \int_{\sigma(T)}^T\sigma^3\|p-\tp\|_4^4 dt.
\ea
\]

\[\ba
\int_0^{\sigma(T)} \sigma\|\g u\|_3^3 dt &\leq C\int_0^{\sigma(T)}\sigma \|\g u\|_2^{\f 32} \|\g u\|_6^{\f 32}dt\\
&\leq C\int_0^{\sigma(T)}\sigma \|\g u\|_2^{\f 32} \big(\|\rho^{\f 12}\|_6\|\rho\dot u\|_{\f 32} + \|p-\tp\|_6\big)^{\f32}dt\\
&\leq C(\mrho)\int_0^{\sigma(T)}\sigma \|\g u\|_2^{\f 32}\|\rho^{\f 12}\dot u\|_2^{\f32}dt + C(\mrho) \int_0^{\sigma(T)}\sigma \|\g u\|_2^{\f 32} dt E_0^{\f 1{4\gamma}}\\
&\leq C(\mrho)\int_0^{\sigma(T)}\big(\sigma^{\f{1-\beta}{2}}\|\g u\|_2\big)\|\g u\|_2^{\f 12} \big(\sigma\|\rho^{\f 12}\dot u\|_2^2\big)^{\f34}dt + C(\mrho) E_0^{\f 34+\f 1{4\gamma}}\\
&\leq C(\mrho, M) E_0^{\f 14} A_1(T)^{\f 34} + C(\mrho) E_0^{\f 34+\f 1{4\gamma}}.
\ea\]
Therefore,
\[\ba
\int_0^T\into \sigma|\g u|^{3} d xdt&\leq C E_0 + C(\mrho )A_2(T)A_1(T) + C(\mrho, M) E_0^{\f 14} A_1(T)^{\f 34} + C(\mrho) E_0^{\f 34+\f 1{4\gamma}}\\
&\qquad + \int_\sigma^T\sigma^3\|p-\tp\|_4^4 dt.
\ea\]
{\bf Step 2.} Estimate of $\int_0^T \sigma^{3}\|p-\tp\|_4^4dt$. Note that
\[
\p_t(p-\tp)+u \cdot \nabla(p-\tp)+\gamma(p-\tp) \div u+\gamma \tp \div u - (\gamma-1)\aint (p-\tp)\div udx=0.
\]
Multiplying the above equality by $3(p-\tp)^2$ yields
\[\ba
\p_t(p-\tp)^3 &+u\cdot\g (p-\tp)^3 +3\gamma (p-\tp)^3\div u +3\gamma\tp(p-\tp)^2\div u\\
&\qquad- 3(\gamma-1)(p-\tp)^2\aint (p-\tp)\div udx = 0.
\ea\]
Integrating on $\O$ shows that
\[\ba
\f{d}{dt}&\into(p-\tp)^3 dx +(3\gamma-1)\into (p-\tp)^3\div u dx +3\gamma\tp\into (p-\tp)^2\div udx\\
&\qquad - 3(\gamma-1)\into(p-\tp)^2dx\aint (p-\tp)\div udx =0
\ea\]
Since
\[
(p-\tp)^3\div u = \f{F + p}{\mu+\ld}(p-\tp)^3 = \f{1}{\mu+\ld} |p-\tp|^4 + \f{F + \tp}{\mu+\ld}(p-\tp)^3.
\]
it holds that
\[\ba
\f{d}{dt}&\into(p-\tp)^3 dx +\f{3\gamma-1}{\mu+\ld}\|p-\tp\|_4^4 \\
&\qquad+ \f{3(\gamma-1)}{\mu+\ld}\into (p-\tp)^3(F+\tp)dx +3\gamma\tp\into (p-\tp)^2\div udx\\
&=3(\gamma-1)\into(p-\tp)^2dx\aint (p-\tp)\div udx\\
%&\leq C\|p-\tp\|_2^2 \Big|\into (p-\tp)\div udx\Big|\\
&\leq C\|p-\tp\|_2^2\|\g u\|_2
\ea\]
Consequently,
\[\ba
\f{3\gamma-1}{\mu+\ld}&\|p-\tp\|_4^4 \leq -\f{d}{dt}\into(p-\tp)^3 dx + C(\mrho)\|p-\tp\|_2^2\|\g u\|_2\\
&\qquad- \f{3(\gamma-1)}{\mu+\ld}\into (p-\tp)^3(F+\tp)dx -3\gamma\tp\into (p-\tp)^2\div udx.
\ea\]
Thus
\[\ba
\f{3\gamma-1}{\mu+\ld} & \sigma^3\|p-\tp\|_4^4 \leq -\f{d}{dt}\Big(\sigma^3\into(p-\tp)^3 dx\Big) + 3\sigma^2\sigma'\into(p-\tp)^3 dx\\
&\qquad- \f{3(\gamma-1)}{\mu+\ld}\sigma^3\into (p-\tp)^3(F+\tp)dx -3\gamma\tp\sigma^3\into (p-\tp)^2\div udx\\
&\qquad + C\sigma^3\|p-\tp\|_2^2\|\g u\|_2.
\ea\]
Integrating above inequality over $[0, T]$ gives
\[
\ba
\f{3\gamma-1}{2(\mu+\ld)} &\int_0^T \sigma^3\|p-\tp\|_4^4 dt = -\sigma^3\into(p-\tp)^3 dx + 3\int_0^{\sigma(T)}\sigma^2\into(p-\tp)^3 dxdt\\
&\qquad- \f{3(\gamma-1)}{\mu+\ld}\int_0^T \sigma^3\into (p-\tp)^3(F+\tp)dxdt \\
&\qquad-3\gamma\tp\int_0^T \sigma^3\into (p-\tp)^2\div udxdt + C\|\g u\|_2^2.
\ea
\]
Therefore,
\[
\ba
\f{3\gamma-1}{2(\mu+\ld)} &\int_0^T \sigma^3\|p-\tp\|_4^4 dt \leq \sup\limits_{0\leq t\leq T}\|p-\tp\|_3^3 + C\int_0^{\sigma(T)}\|p-\tp\|_3^3 dt\\
&\qquad + \eta_1\f{3(\gamma-1)}{\mu+\ld} \int_0^T \sigma^3\|p-\tp\|_4^4dxdt + C_{\eta_1} \f{3(\gamma-1)}{\mu+\ld}\int_0^T \sigma^3\|F+\tp\|_4^4 dt \\
&\qquad + \eta_2 \int_0^T \sigma^3\|p-\tp\|_4^4dxdt + C_{\eta_2}\int_0^T\|\g u\|_2^2 dt.
\ea
\]
Taking $\eta_1=\f 14,\ C_{\eta_1}=2$ and $\eta_2$ small, one gets
\be
\int_0^T \sigma^3\|p-\tp\|_4^4 dt \leq C(\mrho)E^{\f 1\gamma} + \int_0^T \sigma^3\|F+\tp\|_4^4 dt + C E_0.
\ee

Now, by using equation\eqref{ns-0} and \eqref{F-1}, one can get
\[
F + \tp = \f{2\mu}{2\mu+\ld}(A(p)+\tp) + R(\rho\dot u) + B(F) + \f{\ld}{2\mu+\ld}\tp,
\]
\[\ba
\|F+\tp\|_4^4 &\leq 8\big(\f{2\mu}{2\mu+\ld}\big)^4\|A(p)+\tp\|_4^4 + C\|R(\rho\dot u) + B(F) +\f{\ld}{2\mu+\ld}\tp\|_4^4\\
&\leq 8\tC(4)\big(\f{2\mu}{2\mu+\ld}\big)^4\|p-\tp\|_4^4 + C\|R(\rho\dot u) + B(F) +\f{\ld}{2\mu+\ld}\tp\|_4^4
\ea\]
Since $\f{2\mu}{2\mu+\ld}\leq \f 1{4 \tC^{\f 14}(4)}$, it holds that
\[
\int_0^T \sigma^3\|p-\tp\|_4^4 dt \leq C(\mrho)E^{\f 1\gamma} + C \int_0^T \sigma^3\|R(\rho\dot u) + B(F) +\f{\ld}{2\mu+\ld}\tp\|_4^4 dt + C E_0.
\]
Since $R(\rho\dot u) + B(F) +\f{\ld}{2\mu+\ld}\tp =\f{2\mu}{2\mu+\ld}(A(p)+\tp)+(\mu+\ld)\div u -p + \tp$, one has
\[\ba
\|R(\rho\dot u) + B(F) +\f{\ld}{2\mu+\ld}\tp\|_2 & \leq C\|\g u\|_2 + C\|p-\tp\|_2\\
&\leq C\big(\|\g u\|_2 + E_0^{\f 1{2\gamma}}\big)
\ea\]
and by lemma \ref{lem-Ag}
\[
\|\g R(\rho\dot u)\|_{\f 32} \leq C\|\rho\dot u\|_{\f 32}\leq \|\rho^{\f 12}\dot u\|_2\|\rho^{\f 12}\|_6\leq C(\mrho)\|\rho^{\f 12}\dot u\|_2.
\]
Noting that
\[
\aint R(\rho\dot u) + B(F) +\f{\ld}{2\mu+\ld}\tp dx= \aint F + \tp - \f{2\mu}{2\mu+\ld}(A(p)+\tp) dx=0,
\]
one gets
\[\ba
\|R(\rho\dot u) &+ B(F) +\f{\ld}{2\mu+\ld}\tp\|_4^4 \\
&\leq \|R(\rho\dot u) + B(F) +\f{\ld}{2\mu+\ld}\tp\|_2 \|R(\rho\dot u) + B(F) +\f{\ld}{2\mu+\ld}\tp\|_6^3\\
&\leq \|R(\rho\dot u) + B(F) +\f{\ld}{2\mu+\ld}\tp\|_2 \|\g R(\rho\dot u)\|_{\f 32}^3\\
&\leq  C(\mrho)\big(\|\g u\|_2 + E_0^{\f 1{2\gamma}}\big)\|\rho^{\f 12}\dot u\|_2^3
\ea\]
Then,
\[\ba
\int_0^T &\sigma^3 \|R(\rho\dot u) + B(F) +\f{\ld}{\mu+\ld/2}\tp\|_4^4dt
    \leq C(\mrho)\int_0^T \sigma^3\big(\|\g u\|_2 + E_0^{\f 1{2\gamma}}\big)\|\rho^{\f 12}\dot u\|_2^3 dt
\ea\]
\[\ba
\int_0^T &\sigma^3\|\g u\|_2\|\rho^{\f 12}\dot u\|_2^3 dt
    \leq \int_0^T \sigma^{\f 12}\|\g u\|_2 \big(\sigma^{\f 32}\|\rho^{\f 12}\dot u\|_2\big)\big(\sigma\|\rho^{\f 12}\dot u\|_2^2\big) dt\\
&\leq A_1^{\f 32}(T) A_2(T)
\ea\]

\[\ba
\int_0^T& \sigma^3 \|\rho^{\f 12}\dot u\|_2^3 dt
    \leq \int_0^T \big(\sigma^{\f 32}\|\rho^{\f 12}\dot u\|_2\big)\big(\sigma\|\rho^{\f 12}\dot u\|_2^2\big) dt\\
&\leq A_1(T) A_2(T)
\ea\]
Therefore, it holds that
\be
\int_0^T \sigma^{3}\|p-\tp\|_4^4dt\leq
C(\mrho)E^{\f 1\gamma} + C(\mrho) A_1^{\f 32}(T) A_2(T) + C(\mrho)E_0^{\f 1{2\gamma}}A_1(T) A_2(T).
\ee
{\bf Step 3.} Collecting all the estimates above, one conclude that
\[\ba
\int_0^T\into \sigma|\g u|^{3} d xdt&\leq C E_0 + C(\mrho )A_2(T)A_1(T) + C(\mrho, M) E_0^{\f 14} A_1(T)^{\f 34} + C(\mrho) E_0^{\f 34+\f 1{4\gamma}}\\
&\qquad + \int_\sigma^T\sigma^3\|p-\tp\|_4^4 dt.
\ea\]
\be
\int_0^T \sigma^{3}\|p-\tp\|_4^4dt\leq
C(\mrho)E^{\f 1\gamma} + C(\mrho) A_1^{\f 32}(T) A_2(T) + C(\mrho)E_0^{\f 1{2\gamma}}A_1(T) A_2(T).
\ee
\[\ba
A_1(T) + A_2(T)&\leq C E_0^{\f{1}{\gamma}} + C \int_0^T\into \sigma|\g u|^{3} d xdt+ C\int_0^T \sigma^{3}\|p-\tp\|_4^4dt\\
&\leq C E_0^{\f{1}{\gamma}} + C E_0 + C(\mrho )A_2(T)A_1(T) + C(\mrho, M) E_0^{\f 14} A_1(T)^{\f 34} + C(\mrho) E_0^{\f 34+\f 1{4\gamma}}\\
&\qquad + C(\mrho)E^{\f 1\gamma} + C(\mrho) A_1^{\f 32}(T) A_2(T) + C(\mrho)E_0^{\f 1{2\gamma}}A_1(T) A_2(T).
\ea\]
Due to $1<\gamma<2$, one obtains
\[\ba
A_1(T) + A_2(T)&\leq C E_0^{\f{1}{\gamma}} + C E_0 + C(\mrho )A_2(T)A_1(T) \\
&\qquad + C(\mrho, M) E_0^{\f 14} A_1(T)^{\f 34} + C(\mrho) E_0^{\f 34+\f 1{4\gamma}}\\
&\qquad + C(\mrho)E^{\f 1\gamma} + C(\mrho) A_1^{\f 32}(T) A_2(T) + C(\mrho)E_0^{\f 1{2\gamma}}A_1(T) A_2(T)\\
&\leq C E_0^{\f{1}{\gamma}}
+ C(\mrho, M) E_0^{\f 14 + \f 3{8\gamma}} + C(\mrho) E_0^{\f 34+\f 1{4\gamma}}\\
&\leq C_2(\mrho, M)E_0^{\f 1{2\gamma}+\f 1 8}.
\ea\]
Therefore, \eqref{A12est} follows as long as $\eps_4\leq \min\left\{\eps_1, \eps_2, \eps_3, \f{1}{C_2^8(\mrho, M)}\right\}$.
\endproof

%%%%%%%%%%%%%%%%%%%%%%%%%%%%%%%%%%%%%%%%%%%%%%%%%%%%%%%

\section{ Estimate of $\|\rho\|_\q0$}
We start with some weighted estimate for $A(p)$ in the first part of the decomposition.
\begin{lemma}\label{lem-rhoA}
For any $1<\a<\infty$, it holds that
\be\ba
\into|\rho|^\a& |A(p)|dx\leq \tC(\f{\a+\gamma}{\gamma}) \|\rho\|_{\a+\gamma}^{\a+\gamma}.
\ea\ee
\end{lemma}

\proof
\[\ba
\into|\rho|^\a& |A(p)|dx\leq \|\rho^\a\|_{\f{\a+\gamma}{\a}}\|A(p)\|_{\f{\a+\gamma}{\gamma}}
    \leq \tC(\f{\a+\gamma}{\gamma}) \|\rho^\a\|_{\f{\a+\gamma}{\a}} \|p\|_{\f{\a+\gamma}{\gamma}}.
\ea\]
\endproof

The following estimate on the boundary integral of the effective viscous flux F is the key element to derive the higher regularity of the density.
\begin{lemma}\label{Lem-BF}
For any $1<\a<\infty$, it holds that
\be\ba
\f d{dt}&\Big(\into \rho^\a dx\into \rho u\cdot x dx \Big)
%    \leq \int_{\p\O} F(y) dy\into \rho^\a dx + C E_0^{\f 12}\f d{dt} \into \rho^\a dx\\
%&\quad + C\|\g u\|_2^2\into \rho^\a dx +2\into\rho^{\a+\gamma}dx\\
    \leq \int_{\p\O} F(y) dy\into \rho^\a dx\\
    &\qquad + C(\mrho) E_0^{\f 12}(\a-1)\into \rho^\a dx \|\rho^{\f 12}\dot u\|_{2}^{\f 12}\|\g\dot u\|_2^{\f 12}\\
    &\qquad + C\|\g u\|_2^2\into \rho^\a dx +\big(2+C E_0^{\f 12}(\a-1)\big)\into\rho^{\a+\gamma}dx.
\ea\ee
\end{lemma}
\proof Rewrite the momentum balance equation in \eqref{ns-0}
\[\ba
(\rho u)_t + \div(\rho u\otimes u) = \g F -\mu\g\times(\g\times u).
\ea\]
Multiplying $x$ and integrating on $\O$, one gets
\[\ba
\f d{dt}&\into \rho u\cdot x dx + \into \div(\rho u\otimes u)\cdot x = \into x\cdot\g F dx -\mu\into\g\times(\g\times u)\cdot x dx\\
&=\int_{\p\O}F x\cdot n ds - \into 2F dx = \int_{\p\O}F x\cdot n ds +\into 2p dx.
\ea\]
Therefore,
\be\label{IBF}
\ba
\int_{\p\O}F ds &= \int_{\p\O} F x\cdot n ds = \f d{dt}\into \rho u\cdot x dx + \into \div(\rho u\otimes u)\cdot x - \into 2p dx\\
&=\f d{dt}\into \rho u\cdot x dx - \into \rho |u|^2+2pdx.
\ea\ee
Then
\[\ba
\int_{\p\O}& F(y) dy\into \rho^\a dx =\into \rho^\a dx\Big(\f d{dt}\into \rho u\cdot x dx \Big) -\into \rho^\a dx \into \rho |u|^2+2pdx\\
&=\f d{dt}\Big(\into \rho^\a dx\into \rho u\cdot x dx \Big) - \into \rho u\cdot x dx\Big(\f d{dt} \into \rho^\a dx\Big) -\into \rho^\a dx \into \rho |u|^2+2pdx\\
&\geq \f d{dt}\Big(\into \rho^\a dx\into \rho u\cdot x dx \Big) - \into \rho u\cdot x dx\Big(\f d{dt} \into \rho^\a dx\Big) - C\|\g u\|_2^2\into \rho^\a dx -2\into\rho^{\a+\gamma}dx.
\ea\]

On the other hand, direct estimates give
\[\ba
\Big|\into &\rho u\cdot x dx\Big(\f d{dt} \into \rho^\a dx\Big)\Big| \leq C E_0^{\f 12}(\a-1) \into |\rho^\a\div u|dx\\
&\leq C E_0^{\f 12}(\a-1) \into |\rho^\a\div v|+ |\rho^\a\div w|dx\\
&= I + II.
\ea\]
$I$ can be estimated as follows:
\[\ba
I &\leq C E_0^{\f 12}(\a-1)\|\rho^\a\|_{\f{\gamma+\a}{\a}} \|\div v\|_{\f{\gamma+\a}{\gamma}}
\leq C E_0^{\f 12}(\a-1)\|\rho\|_{\a+\gamma}^{\a+\gamma}.
\ea\]
To estimate $II$, one notes that
\[\ba
\into |\rho^\a\div w|dx &\leq \into \rho^\a dx \|\g w\|_\infty
\leq \into \rho^\a dx \|\g w\|_4^{\f 12}\|\g^2 w\|_4^{\f 12}\\
&\leq \into \rho^\a dx \|\g^2 w\|_{\f 43}^{\f 12}\|\g^2 w\|_4^{\f 12}
\leq \into \rho^\a dx \|\rho^{\f 12}\dot u\|_{2}^{\f 12}\|\rho^{\f 12}\|_4^{\f 12}\|\rho\dot u\|_4^{\f 12}\\
&\leq C(\mrho)\into \rho^\a dx \|\rho^{\f 12}\dot u\|_{2}^{\f 12}\|\g\dot u\|_2^{\f 12}
\ea\]
Therefore
\[
\Big|\into \rho u\cdot x dx\Big(\f d{dt} \into \rho^\a dx\Big)\Big| \leq C E_0^{\f 12}(\a-1)\|\rho\|_{\a+\gamma}^{\a+\gamma}
+ C(\mrho) E_0^{\f 12}(\a-1)\into \rho^\a dx \|\rho^{\f 12}\dot u\|_{2}^{\f 12}\|\g\dot u\|_2^{\f 12}.
\]
Consequently, the lemma is proved.
\endproof.

Now we are ready to derive the desired the $L^{q_0}$ bound of the density which closes the bootstrap argument in the key Proposition \ref{prop-1}.
\begin{lemma} \label{lem-rq0}
It holds that
\be
\sup\limits_{0\leq t\leq T} \|\rho\|_{q_0} \leq \f 74\mrho.
\ee
provided that $\f{2\mu}{2\mu+\ld} \leq \f 1{4\tC(\f{\gamma+11}{\gamma-1})}$ and $E_0\leq \eps_5$
with $\eps_5$ is small enough depending on $\mu$, $\ld$, $\gamma$, $q_0$, $\trho$, $\mrho$, $\beta$, $\d_0$, $M$.
\end{lemma}
\def\a{{q_0}}
\proof Recall $\a=\f{12\gamma}{\gamma-1}$, one has
\[
\f{d}{dt}\into \rho^{\a}dx + (\a-1)\into\rho^\a\div u = 0.
\]
Due to
\[
(\mu+\ld)\div u  = F+p = p + \f{2\mu}{2\mu+\ld} A(p) + R(\rho\dot u) + B(F),
\]
one has
\[
\f{d}{dt}\into \rho^{\a}dx + \f{\a-1}{\mu+\ld}\Big[\into\rho^\a p dx+ \f{2\mu}{2\mu+\ld} \into\rho^\a A(p)dx + \into\rho^\a R(\rho\dot u)dx + \into\rho^\a B(F)dx\Big]= 0
\]
Thus, it holds that
\[\ba
\f{d}{dt}&\into \rho^{\a}dx + \f{\a-1}{\mu+\ld}\into\rho^{\a+\gamma} dx + \f{\a-1}{\mu+\ld}\into\rho^\a B(F)dx\\
    &\leq \f{2\mu}{2\mu+\ld}\f{\a-1}{\mu+\ld} \into\rho^\a |A(p)|dx + \f{\a-1}{\mu+\ld}\into\rho^\a |R(\rho\dot u)|dx.
\ea\]
Then, one can obtain that
\[\ba
\f{d}{dt}&\into \rho^{\a}dx + \f{\a-1}{\mu+\ld}\into\rho^{\a+\gamma} dx + \f{\ld}{2\pi (2\mu + \ld)}\f{\a-1}{\mu+\ld}\into\rho^\a dx\int_{\p\O}Fds\\
    &\leq \f{2\mu}{2\mu+\ld}\f{\a-1}{\mu+\ld} \into\rho^\a |A(p)|dx + \f{\a-1}{\mu+\ld}\into\rho^\a |R(\rho\dot u)|dx.
\ea\]
\be\ba
\f d{dt}&\Big(\into \rho^\a dx\into \rho u\cdot x dx \Big)
    \leq \int_{\p\O} F(y) dy\into \rho^\a dx\\
    &\qquad + C(\mrho) E_0^{\f 12}(\a-1)\into \rho^\a dx \|\rho^{\f 12}\dot u\|_{2}^{\f 12}\|\g\dot u\|_2^{\f 12}\\
    &\qquad + C\|\g u\|_2^2\into \rho^\a dx +\big(2+C E_0^{\f 12}(\a-1)\big)\into\rho^{\a+\gamma}dx.
\ea\ee
It follows from lemma \ref{Lem-BF}, one gets
\[\ba
\f{d}{dt}&\into \rho^{\a}dx + \f{\a-1}{\mu+\ld}\into\rho^{\a+\gamma} dx
    + \f{\ld}{2\pi (2\mu + \ld)}\f{\a-1}{\mu+\ld}\f d{dt}\Big(\into \rho^\a dx\into \rho u\cdot x dx \Big)\\
&\leq \f{2\mu}{2\mu+\ld}\f{\a-1}{\mu+\ld} \into\rho^\a |A(p)|dx + \f{\a-1}{\mu+\ld}\into\rho^\a |R(\rho\dot u)|dx\\
    &\qquad + C(\mrho) \f{\ld}{2\pi (2\mu + \ld)}\f{(\a-1)^2}{\mu+\ld}E_0^{\f 12}\|\rho^{\f 12}\dot u\|_{2}^{\f 12}\|\g\dot u\|_2^{\f 12}\into \rho^\a dx \\
    &\qquad + C\f{\ld}{2\pi (2\mu + \ld)}\f{\a-1}{\mu+\ld}\|\g u\|_2^2\into \rho^\a dx +\f{\ld}{2\pi (2\mu + \ld)}\f{\a-1}{\mu+\ld}\big(2+C E_0^{\f 12}(\a-1)\big)\into\rho^{\a+\gamma}dx.
\ea\]
Then, using lemma \ref{lem-rhoA} that
\be\ba
\f{d}{dt}&\into \rho^{\a}dx + \f{\a-1}{\mu+\ld}\into\rho^{\a+\gamma} dx
    + \f{\ld}{2\pi (2\mu + \ld)}\f{\a-1}{\mu+\ld}\f d{dt}\Big(\into \rho^\a dx\into \rho u\cdot x dx \Big)\\
&\leq \f{\a-1}{\mu+\ld}\tC(\f{\a+\gamma}{\gamma})\f{2\mu}{2\mu+\ld} \into\rho^{\a+\gamma}dx + \f{\a-1}{\mu+\ld}\into\rho^\a |R(\rho\dot u)|dx\\
    &\qquad + C(\mrho)\Big( E_0^{\f 12} \|\rho^{\f 12}\dot u\|_{2}^{\f 12}\|\g\dot u\|_2^{\f 12} + \|\g u\|_2^2\Big)\|\rho\|_\a^\a\\
    &\qquad + \f{\ld\big(2+C_1 E_0^{\f 12}(\a-1)\big)}{2\pi(2\mu + \ld)}\f{\a-1}{\mu + \ld}\into\rho^{\a+\gamma}dx.
\ea\ee
By taking $\tC(\f{\a+\gamma}{\gamma})\f{2\mu}{2\mu+\ld} =\tC(\f{\gamma+11}{\gamma-1})\f{2\mu}{2\mu+\ld}\leq \f 14$ and $C_1 \eps_5^{\f 12}(\a-1)\leq 1$, one gets
\[\ba
\f {d}{dt}\|\rho\|^\a_\a &+C_2\f d{dt}\Big(\into \rho^\a dx\into \rho u\cdot x dx \Big)+ \f 14\f{\a-1}{\mu+\ld}\|\rho\|^{\a+\gamma}_{\a+\gamma}\\
    &\leq \f{\a-1}{\mu+\ld}\into\rho^\a |R(\rho\dot u)|dx
    +C(\mrho)\Big( E_0^{\f 12} \|\rho^{\f 12}\dot u\|_{2}^{\f 12}\|\g\dot u\|_2^{\f 12} + \|\g u\|_2^2\Big)\|\rho\|_\a^\a.
\ea\]

Let
\be
y(t) = \|\rho\|^\a_\a + C_2\into \rho^\a dx\into \rho u\cdot x dx.
\ee
Then, by taking $C_2\eps_5^{\f 12}< \f 14$, one can obtain that
\be\label{rhor-1}
\f 34 \|\rho\|^\a_\a \leq y(t)\leq \f 54\|\rho\|^\a_\a,
\ee
and
\be\label{rhor-2}
\ba
y'(t) &+ \f 14\f{\a-1}{\mu+\ld}\|\rho\|^{\a+\gamma}_{\a+\gamma}
\leq \f{\a-1}{\mu+\ld}\into\rho^\a |R(\rho\dot u)|dx\\
    &\qquad +C(\mrho)\Big( E_0^{\f 12} \|\rho^{\f 12}\dot u\|_{2}^{\f 12}\|\g\dot u\|_2^{\f 12} + \|\g u\|_2^2\Big) y(t).
\ea\ee

Next, the following variation of Zlotnik  inequality will be used to get the uniform (in time) upper bound of the density $\rho$ from \eqref{rhor-1}-\eqref{rhor-2}.
\begin{lemma}[\cite{zl1,Huang-Li-Xin}]\label{Lem-zlot}
Let the function $y$ satisfy
\be
 y'(t)\le g(y)+b'(t) \mbox{  on  } [0,T] ,\quad y(0)=y^0,
\ee
with $ g\in C(R)$ and $y,b\in W^{1,1}(0,T).$ If $g(\infty)=-\infty$
and \be\label{a100} b(t_2) -b(t_1) \le N_0 +N_1(t_2-t_1)\ee for all
$0\le t_1<t_2\le T$
  with some $N_0\ge 0$ and $N_1\ge 0,$ then
\be
 y(t)\le \max\left\{y^0,\overline{\zeta} \right\}+N_0<\infty
\mbox{ on
 } [0,T],
\ee
  where $\overline{\zeta} $ is a constant such
that \be\label{a101} g(\zeta)\le -N_1 \quad\mbox{ for }\quad \zeta\ge \overline{\zeta}.\ee
\end{lemma}

\begin{remark}
  In \cite{zl1}, the result holds for equation
  \be
  y'(t)= g(y)+b'(t).
  \ee
  However, following the procedure of his proof, it can be extended to inequality without any difficulty, see \cite{Huang-Li-Xin}.
\end{remark}

Now we are ready to derive uniform $L^q$ estimates for the density.
First, it follows from \eqref{rhor-1}-\eqref{rhor-2} that

\be\label{zlot-1}
\ba
y'(t) &+ \f 14\f{\a-1}{\mu+\ld}c_0y(t)^{\frac{\a+\gamma}{\a}}
  \leq \f{\a-1}{\mu+\ld}\into\rho^\a |R(\rho\dot u)|dx\\
    &\qquad +C(\mrho)\Big( E_0^{\f 12} \|\rho^{\f 12}\dot u\|_{2}^{\f 12}\|\g\dot u\|_2^{\f 12} + \|\g u\|_2^2\Big)y(t)\\
& \le \f{4(\a-1)}{3(\mu+\ld)} \|R(\rho\dot u)\|_{\infty}y(t)
    +C(\mrho)\Big( E_0^{\f 12} \|\rho^{\f 12}\dot u\|_{2}^{\f 12}\|\g\dot u\|_2^{\f 12} + \|\g u\|_2^2\Big) y(t),
\ea
\ee
where one has used the following Holder's inequality:
\be\label{zlot-2}
\into\rho^{\a}dx \le \Big(\into \rho^{\a+\gamma}dx\Big)^{\frac{\a+\gamma}{\a}}|\O|^{\frac{\a+\gamma}{\gamma}}
\defeq c_0\Big(\into\rho^{\a+\gamma}dx\Big)^{\frac{\a+\gamma}{\a}}.
\ee

Hence, then (\ref{zlot-1}) can be rewritten as
\be
D_ty \le f(y) + b'(t),
\ee
where
\be
\ba
& f(y) = - \f 14\f{\a-1}{\mu+\ld}c_0y(t)^{\frac{\a+\gamma}{\a}},\\
&
b(t) = \int_0^t\Big[\f{4(\a-1)}{3(\mu+\ld)} \|R(\rho\dot u)\|_{\infty} + C(\mrho)\Big( E_0^{\f 12} \|\rho^{\f 12}\dot u\|_{2}^{\f 12}\|\g\dot u\|_2^{\f 12} + \|\g u\|_2^2\Big)\Big]y(s)ds.
\ea
\ee

For $t\in[0,\sigma(T)]$, one deduce that for all $0\le t_1<t_2\le\sigma(T)$,
\be
\ba
& |b(t_2)-b(t_1)|\\
& \le C(\mrho)\int_0^{\sigma(T)}\|\rho^{\f 12}\dot u\|_{2}^{\f 12} \|\g \dot u\|_{2}^{\f 12} dt
  + \|\nabla u\|_{2}^2dt\\
& \defeq B_1 + B_2
\ea
\ee
$B_i$ can be estimated as follows:
\be
\ba
B_1 &= \int_0^{\sigma(T)} \|\rho^{\f 12}\dot u\|_{2}^{\f 12} \|\g \dot u\|_{2}^{\f 12} dt
    \leq \int_0^{\sigma(T)} \sigma^{-\f{2-\beta}{4}}\|\rho^{\f 12}\dot u\|_{2}^{\f12} \Big(\sigma^{2-\beta}\|\g \dot u\|_{2}^{2}\Big)^\f{1}{4} dt\\
&\leq \Big(\int_0^{\sigma(T)} \sigma^{-\f{2-\beta}{3}}\|\rho^{\f 12}\dot u\|_{2}^{\f23}dt\Big)^{\f 34} \Big(\int_0^{\sigma(T)}\sigma^{2-\beta}\|\g \dot u\|_{2}^{2}dt\Big)^\f{1}{4}\\
&\leq C(\mrho, M)\Big(\int_0^{\sigma(T)} \sigma^{-\f{2-\beta}{3}}\|\rho^{\f 12}\dot u\|_{2}^{\f23}dt\Big)^{\f 34}\\
&\leq C(\mrho, M)\Big(\int_{0}^{\sigma(T)} \sigma^{-\left[(2-\beta)\left(-\delta_{0}+2 / 3\right)+\delta_{0}\right]}
    \Big(\sigma^{2-\beta}\|\rho^{1 / 2} \dot{u}\|_{2}^{2}\Big)^{-\delta_{0}+1 / 3}\Big(\sigma\|\rho^{1 / 2} \dot{u}\|_{2}^{2}\Big)^{\delta_{0}} d t\Big)^{\f 34}\\
&\leq C(\mrho, M)(A_1(\sigma(T)))^{\f {3\d_0}{4}}.
\ea
\ee
and
\[
B_2 = \int_0^T \|\g u\|_2^2 dx \leq E_0^{\f 12}.
\]
Thus, it holds that
\be
|b(t_2)-b(t_1)| \le C(\mrho, M)(A_1(\sigma(T)))^{\f {3\d_0}{4}} + CE_0^{\f 12}
\ee
provided that $E_0\le \eps_5$. Therefore, one can choose $N_0$ and $N_1$ as follows:
\be
N_1 =0, N_0 = C(\mrho, M)(A_1(\sigma(T)))^{\f {3\d_0}{4}} + CE_0^{\f 12}, \bar{\eta} = 1.
\ee
Then
\be
f(\eta) =  - \f 14\f{\a-1}{\mu+\ld}c_0\eta^{\frac{\a+\gamma}{\a}}\le -N_1 =0 \qquad\mbox{for all $\eta\ge\bar{\eta}=1$.}
\ee
Thus lemma \ref{Lem-zlot} implies that
\be
y(t)\le \max\{\f 54\mrho^{\a},1\} + N_0 \le \f 54\mrho^{\a} + C(\mrho, M)(A_1(\sigma(T)))^{\f {3\d_0}{4}} + CE_0^{\f 12}\le \frac{21}{16}\mrho^{\a}.
\ee
provided that C$E_0^{\f 12} + C(\mrho, M)(A_1(\sigma(T)))^{\f {3\d_0}{4}}\le\frac{1}{16}\mrho^{\a}$.

  On the other hand, for $t\in[\sigma(T),T]$, one can get
\be
\ba
|b(t_2)-b(t_1)|& \le C(\mrho)\int_{\sigma(T)}^{T} \|\rho^{\f 12}\dot u\|_{2}^{\f 12} \|\g \dot u\|_{2}^{\f 12} dt
 + C(\mrho)\int_{\sigma(T)}^{T}\|\nabla u\|_{2}^2dt\\
&\le \f 14\f{\a-1}{\mu+\ld}c_0(t_2-t_1) + C(\mrho)E_0 + C(\mrho)\int_{\sigma(T)}^{T}\|\rho^{\f 12}\dot u\|_{2}^{2} \|\g \dot u\|_{2}^{2} dt\\
&\le \f 14\f{\a-1}{\mu+\ld}c_0(t_2-t_1) + C(\mrho)E_0 + C(\mrho)E_0^{\frac{1}{10\gamma}}\int_{\sigma(T)}^{T}\|\nabla{\dot{u}}\|_{2}^2dt\\
&\le \f 14\f{\a-1}{\mu+\ld}c_0(t_2-t_1) + C(\mrho)E_0^{\frac{1}{5\gamma}}
\ea
\ee
provided that $E_0\le \eps_5$. Therefore, we can choose $N_1$ and $N_0$ as follows:
\be
N_1 = \f 14\f{\a-1}{\mu+\ld}c_0, N_0 = C(\mrho)E_0^{\frac{1}{5\gamma}}.
\ee
Note that
\be
f(\eta) = - \f 14\f{\a-1}{\mu+\ld}c_0\eta^{\frac{\a+\gamma}{\a}}\le -N_1, \qquad\mbox{for all $\eta\ge 1$.}
\ee
Thus, lemma \ref{Lem-zlot} yields that
\be
y(t)\le \max\{\f 54\mrho^{\a},1\} + N_0\le \f 54\mrho^{\a} + C(\bar{\rho})E_0^{\frac{1}{5\gamma}}\le \frac{21}{16}\mrho^{\a},
\ee
provided that
\be
E_0\le \eps_5\leq \Big(\frac{\mrho^\a}{16C(\mrho)}\Big)^{5\gamma}.
\ee
Thus lemma \ref{lem-rq0} is proved.
\endproof

Now Proposition \ref{prop-1} follows from lemma \ref{lem-A123} and lemma \ref{lem-rq0}. Hence Theorem \ref{th-main} follows from Proposition \ref{prop-1} by the standard bootstrap argument.

\section*{Acknowledgements}
 X.-D. Huang is partially supported by CAS Project for Young Scientists in Basic Research, Grant No.YSBR-031 and NNSFC Grant No. 11688101. Xin is partially supported by Zheng Ge Ru Foundation, Hong Kong RGC Earmarked Research Grants CUHK-14301421, CUHK-14301023, CUHK-14300819 and CUHK-14302819, and the key projects of NSFC Gronts No. 12131010 and No. 11931013. W. YAN is partially supported by NNSFC Grant No. 11371064 and No. 11871113.

\end{document}